\documentclass[11pt]{article}
	
	\newcommand{\blind}{0}
	
    \makeatletter
    \renewcommand\section{\@startsection {section}{1}{\z@}%
                                       {-3.5ex \@plus -1ex \@minus -.2ex}%
                                       {2.3ex \@plus.2ex}%
                                       {\normalfont\fontfamily{phv}\fontsize{16}{19}\bfseries}}
    \renewcommand\subsection{\@startsection{subsection}{2}{\z@}%
                                         {-3.25ex\@plus -1ex \@minus -.2ex}%
                                         {1.5ex \@plus .2ex}%
                                         {\normalfont\fontfamily{phv}\fontsize{14}{17}\bfseries}}
    \renewcommand\subsubsection{\@startsection{subsubsection}{3}{\z@}%
                                        {-3.25ex\@plus -1ex \@minus -.2ex}%
                                         {1.5ex \@plus .2ex}%
                                         {\normalfont\normalsize\fontfamily{phv}\fontsize{14}{17}\selectfont}}
    \makeatother
	\usepackage{amsmath}
    \usepackage{calc}          
	\usepackage{graphicx}
	\usepackage{enumerate}
	\usepackage{url} 
	
	\usepackage{amsfonts,amssymb,amsthm}
	\usepackage{hyperref}       
    \usepackage{nicefrac}       
    \usepackage{mathtools}
    \usepackage{xparse}         
    \usepackage{graphicx}
    \usepackage{color}
    \usepackage{bm}

    \allowdisplaybreaks
	
	\usepackage{amsfonts,amssymb,amsthm}
	\usepackage{hyperref}       
    \usepackage{nicefrac}       
    \usepackage{mathtools}
    \usepackage{xparse}         
    \usepackage{graphicx}
    \usepackage{bm}
    \newtheorem{thm}{Theorem}
    \newtheorem{lemma}{Lemma}
    \newtheorem{cor}{Corollary}
    \newtheorem{prop}{Proposition}

    \newcommand{\sS}{(s,S)}
    \newcommand{\sSstar}{(s^*,S^*)}

    \newcommand{\sss}{\scriptscriptstyle}
    \newcommand{\Zb}{\mathbb{Z}}
    \newcommand{\Rb}{\mathbb{R}}
    \newcommand{\cN}{\mathbb{N}}
    \newcommand{\Zp}{\Zb_+}
    \newcommand{\Rp}{\Rb_+}
    \newcommand{\la}{\lambda}
    \newcommand{\tht}{\theta}
    \newcommand{\al}{\alpha}
    \newcommand{\ga}{\gamma}
    \newcommand{\io}{\iota}
    \newcommand{\bt}{\beta}

    \newcommand{\Ga}{\Gamma}
    \newcommand{\Dlt}{\Delta}
    \newcommand{\La}{\Lambda}
    \newcommand{\Tht}{\Theta}

    \newcommand{\bal}{\boldsymbol{\al}}
    \newcommand{\bla}{\boldsymbol{\la}}
    \newcommand{\bmu}{\boldsymbol{\mu}}
    \newcommand{\bth}{\boldsymbol{\tht}}
    \newcommand{\bxi}{\boldsymbol{\xi}}
    \newcommand{\bpi}{\boldsymbol{\pi}}
    \newcommand{\bq}{\boldsymbol{q}}
    \newcommand{\bu}{\boldsymbol{U}}
    \newcommand{\bv}{\boldsymbol{v}}
    \newcommand{\wPhi}{\widetilde{\Phi}}
    \newcommand{\barla}{\Bar{\la}}

    \newcommand{\bp}{\boldsymbol{p}}
    \newcommand{\bZ}{\boldsymbol{0}}
    \newcommand{\mcI}{\mathcal{I}}
    \newcommand{\mcS}{\mathcal{S}}
    \newcommand{\mcL}{\mathcal{L}}
    \newcommand{\Sphi}{S_{\Phi}}
    \newcommand{\Sphiph}{\Sphi^{ph}}
    \newcommand{\prR}[1]{p_{_{#1}}}
    \newcommand{\bprR}[1]{\boldsymbol{p\!}_{_{#1}}}
    \newcommand{\binv}{B_{inv}}

    \newcommand{\diagm}[1]{\Dlt\!\apar{#1}}

    \newcommand{\ablk}[2]{A_{#1}^{^{(#2)}}}
    \newcommand{\tablk}[2]{\tilde{A}_{#1}^{^{(#2)}}}
    \newcommand{\ablkij}[3]{\left[A_{#1}^{^{(#2)}}\right]_{#3}}
    \newcommand{\tablkij}[3]{\left[\tilde{A}_{#1}^{^{(#2)}}\right]_{#3}}
    \newcommand{\apar}[1]{\left(#1\right)}
    \newcommand{\abrk}[1]{\left[#1\right]}

    \newcommand{\bdiagi}{\boldsymbol{d}_{_{\io}}(R)}
    \newcommand{\bdiagb}{\boldsymbol{d}_{_{\bt}}}
    \newcommand{\bdiagg}{\boldsymbol{d}_{_{\ga}}}
    \newcommand{\tA}{\tilde{A}}
    \newcommand{\tM}{\tilde{M}}
    \newcommand{\tS}{\tilde{S}}
    \newcommand{\tD}{\tilde{D}}
    \newcommand{\tTh}{\tilde{\Tht}}
    \newcommand{\tGa}{\tilde{\Ga}}
    \newcommand{\tLa}{\tilde{\La}}
    \newcommand{\tbpi}{\tilde{\bpi}}

    \newcommand{\Ai}[1]{A^{^{(#1)}}}
    \newcommand{\tAi}[1]{\tilde{A}^{^{(#1)}}}
    
    \newcommand{\bpist}{\tbpi^{\ast}}
    \newcommand{\inv}[1]{#1^{^{-1}}}
    
    \newcommand{\trps}[1]{{#1}^{\intercal}}
    \newcommand{\pyblk}[2]{\tilde{\bpi}_{_{#1#2}}}
    \newcommand{\pyblkh}[2]{\tilde{\bpi}_{h_{#1#2}}}

    \DeclareDocumentCommand{\eins}{o}{\boldsymbol{e}\IfNoValueTF{#1}{}{_{_{#1}}}}
    \DeclareDocumentCommand{\diagi}{o}{d_{_{\io}}\!(\IfNoValueTF{#1}{z}{#1},R)}
    \DeclareDocumentCommand{\diagb}{o}{d_{_{\bt}}\!(\IfNoValueTF{#1}{z}{#1})}
    \DeclareDocumentCommand{\diagg}{o}{d_{_{\ga}}\!(\IfNoValueTF{#1}{z}{#1})}
    \DeclareDocumentCommand{\E}{m o}{\mathbb{E}\!
        \left[\IfNoValueTF{#2}{}{\left.}#1\IfNoValueTF{#2}{}{\,\right|#2}\right]}
    \DeclareDocumentCommand{\prb}{m o}{\mathbb{P}\!
        \left(#1\IfNoValueTF{#2}{}{\,\vert\,#2}\right)}
        \DeclareDocumentCommand{\sett}{m o}
        {\left\{#1\IfNoValueTF{#2}{}{\,:\,#2}\right\}}
    \DeclareDocumentCommand{\py}{m o o}
        {\tilde{\IfNoValueTF{#3}{\bpi}{\pi}}_{_{#1\IfNoValueTF{#2}{\bullet}{#2}
        \IfNoValueTF{#3}{\bullet}{#3}}}}
	
\begin{document}
	
		\def\spacingset#1{\renewcommand{\baselinestretch}%
			{#1}\small\normalsize} \spacingset{1}
		
		\if0\blind
		{
			\title{\bf A Markov-Modulated $\sS$ Inventory System 
                   with Repeated Calls and Blocked Demands}
			\author{James Cordeiro$^{a*}$, Ying-Ju Chen$^{b*}$, Andr\'{e}s
			  Larrain-Hubach$^{c*}$, and Mark Abramson$^{d**}$ \\\small
		  $^*$ Department of Mathematics, University of Dayton, Dayton, Ohio USA \\\small
		  $^*$ Department of Mathematics, Utah Valley University, Orem, Utah USA \\\small
             $^a$ \emph{Corresponding Author}, cordeirj@yahoo.com, \\\small
             $^b$ ychen4@udayton.edu \\\small
             $^c$ alarrainhubach1@udayton.edu\\\small
             $^d$ mark.abramson@uvu.edu}
        
			\date{}
			\maketitle
		} \fi
		
		\if1\blind
		{
            \title{\bf An $\sS$ Inventory System with Demand Blocking and
                  Retrials in a Random Environment}
			\author{Author information is purposely removed for double-blind review}
			
            \bigskip
			\bigskip
			\bigskip
			\begin{center}
				{\LARGE\bf An $\sS$ Inventory System with Demand Blocking and
                  Retrials in a Random Environment}
			\end{center}
			\medskip
		} \fi
		\medskip
		
	\begin{abstract}
        In this article, we consider a continuous review $\sS$ 
        inventory system with failures of demand fulfillment (service)
        modeled as a Markov-modulated retrial queueing system. The inventory system
        features a single product that experiences Markovian inter-demand and service
        intervals with random service interruptions and instantaneous
        replenishments. A recently developed criterion for the ergodicity of a class of
        discrete-time level-dependent-quasi-birth-and-death (LDQBD) processes with
        convergent transition matrix rows is applied to
        the jump chain of the process in order to elicit a closed-form
        traffic-intensity formula. An analytic solution for the steady-state average
        minimum cost is provided.
	\end{abstract}
			
	\noindent%

	\spacingset{1.5} 

\section{Introduction}

The classical $\sS$ inventory model, which was first investigated in
Arrow and Harris \cite{arrow1951}, and its later variants were
developed to address the practical concerns of inventory management,
and in doing so, posed interesting theoretical questions about
model stability and optimal control. In particular, as Fisher and Hornstein
\cite{fisher_horn_2000} assert, $\sS$ models were extensively studied for use in
retail applications due to the assumption of fixed ordering costs. The intuitive
operation of $\sS$ models, together with their practical relevance, have
given them prominence in the inventory literature.

In the classical $\sS$ inventory model, single demands for a type of item
arrive to the system, and they are fulfilled as long as the inventory
contains at least one item.
However, if there arises an order that depletes the inventory to a level at or below
a critical threshold value $s>0$, then an order for just enough items
to restore the level of the inventory to its maximum
capacity of $S>0$ is made. A time delay between
replenishment orders and deliveries may be specified whenever the
threshold level $s$ is attained, or none at all, as is the case in
what is
termed an instantaneous replenishment model. Observations of the
product level needed to trigger successive replenishment of the inventory
may take place continuously or otherwise over time. The
first $\sS$ models provided for continuous monitoring of inventory
levels, such as in the model of this paper, hence the designation
continuous review. This form of monitoring is the
one that most often characterizes Markovian queueing inventory models.

The further imposition of a queueing model framework to inventory 
systems allows the modeler to leverage analytical techniques
developed for the performance analysis of queues in steady-state
operation. In a Markovian queueing system,
incoming demands are often represented as a Poisson input stream
and their subsequent processing as the in-service durations that are
associated with one or more servers. In the event of
blocked demands due to failures or busy periods of a server,
the retrial queueing models of Artalejo and Krishnamoorthy \cite{artalejo_etal06}, 
and Ushakumari \cite{ushakumari06} may be employed. In retrial models, blocked demands
are redirected into a holding area called a retrial orbit, upon which
each demand persistently reattempts fulfillment at i.i.d. time intervals.
Demands are thus retained in the system without backlogging, i.e., without
a promise of fulfillment, such as happens when items are back-ordered.
Consequently, fulfillment will occur only when items are available and the
ordering system is functioning, as usually occurs in online ordering
scenarios.

In addition to imperfect service, another feature intended to 
free queueing systems from restrictive simplifying assumptions is the
specification of a fluctuating random environment, which was first studied
by Yechiali and Naor \cite{yech_naor71} and expanded upon by Neuts
\cite{neuts71}. This is an independently evolving exogenous stochastic process that modifies the distributional parameters of the various
time durations at evolutionary epochs. Such queueing systems, which are alternatively
referred to as Markov modulated queueing systems, also appear in the context of
queueing inventory systems,
such as in the publications of Karlin \cite{Karlin1960} and Iglehart and
Karlin \cite{IglehartKarlin1963}. Subsequently, the first to study an $\sS$ inventory system with a compound-Poisson demand process modulated by a finite-state Markovian random environment was Feldman \cite{Feldman1978}. Other inventory models that utilize a random environment include, but are not limited to Song and Zipkin \cite{Song1993},
\:{O}zekici and Parlar \cite{Ozekici1999}, and Perry and Posner \cite{Perry2002}.

A notable outcome of the study of Markov-modulated queueing systems is that
their underlying Markov chains were found to be quasi-birth-and-death (QBD)
processes, which are discrete- or continuous-time Markov chains whose
transition matrix entries in block form are arranged according to a distinctive tri-diagonal pattern, as described by the seminal work of Neuts
\cite{neuts1:78,neuts2:78}, who also gives an analytic criterion for their positive recurrence. However, this criterion is limited to QBDs whose transition matrices possess infinitely repeating block rows, save for a finite number of
boundary rows. Such QBDs are termed homogeneous or level-independent QBDs. These are in turn subsumed within a general class of QBDs whose rows do not
repeat, and which are accordingly termed level-dependent QBDs, or LDQBDs.
Markovian inventory models with underlying LDQBDs may be found in Artalejo et al.
\cite{artalejo_etal06}, Ushakumari \cite{ushakumari06}, Krishnamoorthy, Nair,
and Narayanan \cite{Krishnamoorthy2012}, and Ko \cite{Ko2020}.
Analytical criteria for the ergodicity and non-ergodicity of LDQBDs were eventually discovered by Cordeiro,
Kharoufeh, and Oxley \cite{cordeiro2019} for irreducible processes whose transition matrices exhibit element-wise
row convergence to a single limiting block row, which we shall henceforth term row-convergent LDQBDs. Such behavior characterizes a plethora of useful queueing
models, to include the inventory model that is considered in this paper.

To the best of the authors' knowledge, the ergodicity criteria of Cordeiro et al.
\cite{cordeiro2019} has not yet been utilized to develop criteria for the stability
of queueing inventory models whose underlying Markov chains may be classified as LDQBDs.
Therefore, in this paper, we seek to address this concern by formulating
a general traffic intensity formula application using the matrix analytic
approach of Cordeiro et al. \cite{cordeiro2019}. In addition, a means to evaluate the performance
characteristics of such models in steady-state is likewise developed.

The remainder of this paper is organized as follows. Section \ref{sec:LDQBD} introduces
the LDQBD and the drift criterion for the ergodicity of row-convergent LDQBDs. 
After a description of the instantaneous replenishment $\sS$ inventory system in Section \ref{sec:model},
Section \ref{sec:model_stability} establishes that its underlying LDQBD is row-convergent, upon which
an analytic traffic intensity formula for the system is derived using the method of Cordeiro et al.
\cite{cordeiro2019}. With a means to determine positive recurrent inventory systems in hand,
Section \ref{sec:measures} develops steady-state average performance measures
for positive-recurrent systems. Lastly, in Section \ref{sec:optimization}, a comparison of average
cost solutions of stable systems over systems of varying traffic intensity is presented.

\section{Level-Dependent Quasi-Birth-and-Death Processes} \label{sec:LDQBD}

A continuous-time level-dependent quasi-birth-and-death (LDQBD) process is
a bivariate continuous-time Markov chain (CTMC)
$\Phi=\sett{(X(t),Y(t))}[t\ge 0]$ with state space
\[S_{\Phi}=\sett{(i,j)}[i\in\Zp,\,j\in\sett{1,\dots,K}],\]
where $\Zp$ is the set of non-negative integers and $K<\infty$ is some positive
integer value. The $x$-coordinate of
$S_{\Phi}$ is denoted as the level of the process while the $y$-coordinate is
the phase. The infinitesimal generator $Q^*$ of $\Phi$ consists of $K\times K$
block entries that are arrayed in the distinctive tridiagonal form given by
\begin{equation}\label{eq:Q*}\small
    Q^*=[q^*_{ij}]=\begin{bmatrix}
    \ablk{1}{0} & \ablk{0}{0} & 0         & 0         & 0 & \cdots \\
    \ablk{2}{1} & \ablk{1}{1} & \ablk{0}{1} & 0         & 0 & \cdots \\
    0         & \ablk{2}{2} & \ablk{1}{2} & \ablk{0}{2} & 0 & \cdots \\
    0         & 0         & \ablk{2}{3} & \ablk{1}{3} & \ablk{0}{3}     & \cdots \\
    0         & 0           & 0           & \ablk{2}{4} & \ablk{1}{4} & \cdots \\
    \vdots & \vdots & \vdots & \vdots & \vdots & \ddots \\
  \end{bmatrix}
\end{equation}
where 0 denotes the zero matrix and the nonzero entries $\ablk{k}{i}$ vary according
to the level $i$ for each $k\in\sett{0,1,2}$. If the block entries $\ablk{k}{i}$ are
invariant over all levels, that is, $\ablk{k}{i}=A_k$ for all levels $i$,
save for a finite number of initial levels beginning
with level 0, then the process is termed a level-independent, or homogeneous
QBD. The closed-form ergodicity criterion for an irreducible continuous-time
homogeneous QBD, which was derived by Neuts \cite{neuts1:78}, is that the process
is positive-recurrent if and only if
\begin{equation}\label{eq:homog_drift}
  D=\bpi(A_0-A_2)\eins<0,
\end{equation}
where $\eins$ is a column vector of the appropriate dimension (in this case $m$)
whose scalar entries consist entirely of ones
and $\bpi$ is a $m$-dimensional row vector that solves the linear system
$\bpi(A_0+A_1+A_2)=0$ and $\bpi\eins=1$. In either case, the process is referred
to as skip-free, in deference
to the characteristic that no transition of the process may exceed one
level in either the positive or negative direction.

Next, we consider the discrete-time Markov chain (DTMC)
\[\wPhi=\sett{(X_n,Y_n)}[n\in\Zp]\]
with state space $S_{\Phi}$ that is embedded at transitions of the CTMC
$\Phi$; the transition times are enumerated according to $n\in\Zp$.
This is known as the jump process of $\Phi$. Its transition
probability matrix $\tilde{P}$ exhibits the same tridiagonal block structure
\begin{equation}\label{eq:P}\small
  \tilde{P}=\begin{bmatrix}
    \tablk{1}{0} & \tablk{0}{0} & 0          & 0          & 0 & \cdots \\
    \tablk{2}{1} & \tablk{1}{1} & \tablk{0}{1} & 0          & 0 & \cdots \\
    0          & \tablk{2}{2} & \tablk{1}{2} & \tablk{0}{1} & 0 & \cdots \\
    0          & 0          & \tablk{2}{3} & \tablk{1}{3} & \tablk{0}{3} & \cdots \\
    0          & 0          & 0          & \tablk{2}{4} & \tablk{1}{4} & \cdots \\
    \vdots & \vdots & \vdots & \vdots & \vdots & \ddots
  \end{bmatrix}.
 \end{equation}
The elements of $\tablk{k}{i}$ for $k=0,1,2$ and for each
$(i,j)\in S_{\Phi}$ and $j'\in\sett{1,\dots,K}$ are the probabilities
\begin{align*}
  \tablkij{2}{i}{jj'} &= \prb{X_{n+1}=i-1,\,Y_{n+1}=j'}[X_n=i,\,Y_n=j],\quad i\ge 1, \\\notag
  \tablkij{1}{i}{jj'} &= \prb{X_{n+1}=i,\,Y_{n+1}=j'}[X_n=i,\,Y_n=j],
    \quad\tablkij{1}{i}{jj}=0, \quad j'\neq j, \\
  \tablkij{0}{i}{jj'} &= \prb{X_{n+1}=i+1,\,Y_{n+1}=j'}[X_n=i,\,Y_n=j].
\end{align*}

For the purpose of determining system stability, it is necessary to restrict
our attention to the class of irreducible discrete-time LDQBD processes $\wPhi$
for which the following element-wise limits
\begin{equation}\label{asmp:conv}
  \tA_k^*=\lim_{i\to\infty}\tablk{k}{i},\quad\text{ exist for } k=0,1,2,
\end{equation}
and, in addition,
\begin{equation}\label{asmp:stochA}
  \tA^*=\tA_0^*+\tA_1^*+\tA_2^*\;\text{ is a stochastic matrix.}
\end{equation}
In other words, the rows of transition probability matrix of $\wPhi$, which
is subject to Eqns. \eqref{asmp:conv} and \eqref{asmp:stochA}, approach a
limiting row as the level increases. We henceforth term such discrete-time
QBDs as row-convergent LDQBDs.
As described in Cordeiro et al. \cite{cordeiro2019}, the discrete-time row-convergent
LDQBD $\wPhi$ is positive-recurrent if and only if
\begin{equation}\label{eq:ldqbd_erg}
  \tD^*<0,
\end{equation}
where we define the average drift $\tilde{D}^*$ of process $\wPhi$ to be the scalar
quantity
\begin{equation}\label{eq:D*}
 \tD^*=\tbpi^*\apar{\tA_0^*-\tA_2^*}\eins
\end{equation}
and $\tbpi^*$ is the unique $K$-dimensional vector that solves the linear system
\begin{equation}\label{eq:pi*_system}
  \tbpi^*\tA^*=\tbpi^*,\qquad \tbpi^*\cdot\eins=1, \qquad
    \tA^*=\tA^*_0+\tA^*_1+\tA^*_2.
\end{equation}

\section{Model Description}\label{sec:model}

\begin{figure}[!http]
    \centering
    \includegraphics[width=\textwidth]{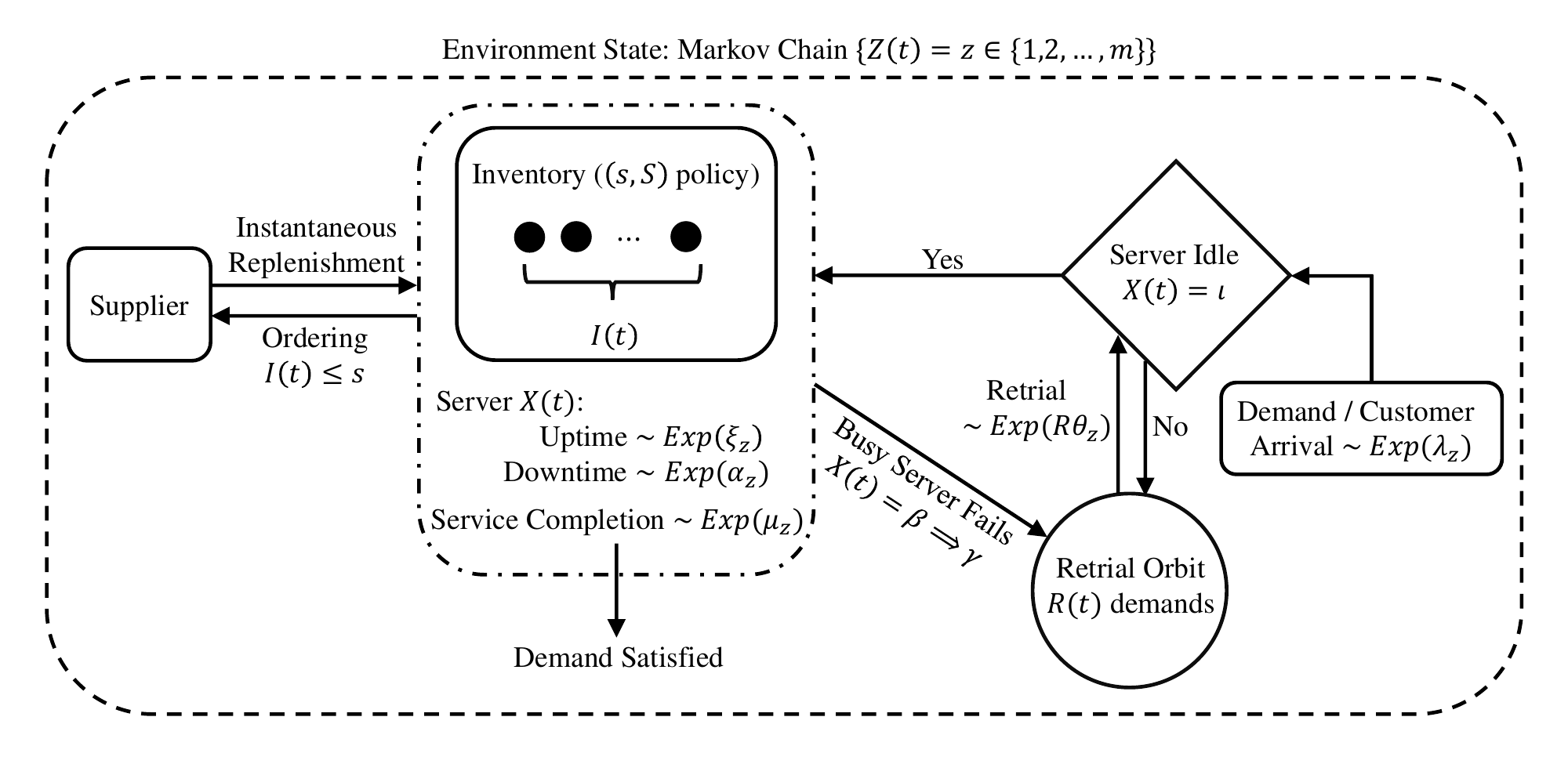}
    \caption{Model Illustration}
    \label{fig:model}
\end{figure}

The system that we consider here (refer to Figure \ref{fig:model})
is a continuous-review $\sS$ inventory system that consists of a single-product
storage facility and a single server that processes incoming demands.
Letting $\cN=\sett{1,2,3,\dots}$, we define $S\in\cN$
to be the fixed inventory storage capacity and $0\le s<S$ to be the threshold
level at which a replenishment of the inventory is triggered. If the level of
product in the inventory drops to the threshold level of $s$, an instantaneous
replenishment of $S-s$ items occurs. Such a replenishment policy maintains the
inventory level in the range $[s+1,S]$, which enforces the requirement that
only one replenishment takes place at any instant of time.

In order to consider the mathematical performance measures
of the system in equilibrium, we will model this inventory system as a
standard $M/M/1$ retrial queueing system with a Poisson arrival
stream of demands that possesses an average interarrival duration
of $1/\la$ and a single server that processes incoming demands according
to exponential service durations that are of average length $1/\mu$. In
lieu of a standard FIFO queue, there is, instead, a retrial orbit with
unrestricted capacity for unsatisfied demands that proceed here from
a busy or failed server. While in orbit, each
of these demands will re-attempt service independently of all other
demands in orbit at
intervals distributed exponentially with average length $1/\tht$. This
results in a combined output stream with an inter-retrial duration that
is distributed exponentially, with an average
duration of $1/(R\tht)$, where $R$ is the current number in orbit.

Before any incoming demand is satisfied, it must be processed by the system
server. The server is assumed at all times to be in one of three states,
namely idle and operational, busy and operational, or failed.
A server that is failed will not satisfy a demand for the product.
The server remains operational for an exponential duration with an average
length $1/\xi$, after which it is considered to be in a failed state.
Repair of the server commences immediately for an exponential duration of the average
length $1/\al$, after which the server is returned to a fully operational and idle
state.

At the time $t=0$, it is assumed that the server is idle and
operational, the inventory is at its maximum level $I(0)=S$, and there
are no demands in the system. Thereafter, single demands arrive to the server
according to the specified Poisson process. If the server
is idle, processing of the incoming demand commences, and the server assumes
a busy state. If the server does not fail, then
the inventory is decremented by one unit at the end of the service duration and the demand then leaves the system.
Subsequently, if the inventory decrements to $s$ items, then an instantaneous
replenishment to the full capacity $S$ of the inventory takes place.

On the other hand, if a demand encounters a busy or a failed server,
it will proceed directly to the retrial orbit.
Likewise, if the server fails while in a busy state, the demand
being processed will immediately proceed to the retrial orbit. In
either case, the number in inventory will not be decremented.
A demand in orbit may obtain service only when the combined retrial
duration with rate $R\theta$ ends when the server is idle. Afterward,
the orbit size is decremented to $R-1$ and a busy period of the server
commences.

We seek to emulate the effect of external influences, such as
fluctuations in economic conditions, by the inclusion of a random environment that varies the exponential distributions of inter-demand arrival times, its subsequent processing (service) times, service up- and down-times,
and times between retrials of service. Accordingly,
the random environment process will be defined here as a finite-state
irreducible CTMC $\sett{Z(t)}[t\in\Rp]$ with state space
$\mcS=\sett{1,\dots,m},m\in\cN$, and
infinitesimal generator $Q=[q_{zz'}]_{z,z'\in\mcS}$. In the standard way,
we denote the total rate out of state $z\in \mathcal{S}$ as
$$q_z=-q_{zz}=\sum_{z'\ne z}q_{zz'},\quad z\in\mcS.$$
If $Z(t)=z$ at a time instant $t\ge 0$, then the exponential parameters
of each process appear as follows:
\begin{center}\small\begin{tabular}{|l|*{6}{c}|}\hline
    \textbf{Process} & Arrival & Service & Uptime & Downtime & Retrial & Environment\\\hline
    \textbf{Rate} & $\la_z$ & $\mu_z$ & $\xi_z$ & $\al_z$ & $\tht_z$ & $q_z$ \\\hline
\end{tabular}\end{center}
For convenience, the parameters are expressed as entries of the respective
$m$-vectors $\bla$, $\bmu$, $\bxi$, $\bal$, $\bth$, and $\bq$.

We next define the random variables that reflect the state of the system
at time $t\ge 0$. Let\vspace{-1ex}
\begin{align*}
  R(t) &= \text{the number of demands in orbit at time } t, \\
  I(t) &= \text{the number of items in the inventory at time } t, \\
  X(t) &= \text{the status of the server at time $t$}, \\
       &=
        \begin{cases}
           \io & \text{ if the server is idle}, \\
           \bt & \text{ if the server is busy}, \\
           \ga & \text{ if the server is failed},
        \end{cases} \\
  Z(t) &= \text{the state of the random environment at time } t.
\end{align*}
Due to the fact that all of the time durations of the process are exponentially distributed, the Markov property holds. Consequently, we may define the system as the
Markov chain
\[\Phi=\sett{(R(t),I(t),X(t),Z(t))}[t\in\Rp]\]
with the state space
\[\Sphi=\sett{(R,I,X,Z)}[R\in\Zp,\,I\in[s+1,S]\cap\Zp,\,
  X\in\{\io,\bt,\ga\},Z\in\mcS].\]
For convenience, we define the finite phase state partition of $\Sphi$
as the set
\[\Sphiph=\sett{(I,X,Z)}[\,I\in[s+1,S]\cap\Zp,\,
  X\in\{\io,\bt,\ga\},Z\in\mcS].\]
If the elements of this set are enumerated in lexicographic order as
\[\mcL=\sett{1,2,\dots,K},\quad K=3m(S-s)\]
we may then rewrite the state space as
\[\Sphi=\sett{(R,k)}[R\in\Zp,\,k\in\mcL].\]
Moreover, the process $\Phi$ possesses an infinitesimal generator matrix
$Q^*=[q^*_{yy'}]$, where both $y=(R,I,X,Z)$ and
$y'=(R',I',X',Z')$ belong to $\Sphi$.
The rows and columns of the matrix are arranged according in the
lexicographic order of ascending orbit size $R$ (level) and the
order given in $\mcL$ at each level $R$. The matrix consequently
appears as in Eqn. \eqref{eq:Q*}.

Elements of the generator matrix $Q^*$ of $\Phi$ will next be specified.
For the purpose of simplification, we define for each $R\in\Zp$ the scalar
values
\begin{align*}
  \diagi &= q_z+\la_z+\xi_z+R\tht_z \\
  \diagb &= q_z+\la_z+\mu_z+\xi_z \\
  \diagg &= q_z+\la_z+\al_z.
\end{align*}
The resulting entries of $\ablk{k}{R}$ for $k=0,1,2$ and
at each level $R$ are depicted in Table \ref{tab:A_2^(R)}.
\begin{table}[tbh]
\begin{center}\small
\begin{tabular}{|c|c|l|l|c|l|}\hline
    $k$ & $R$ & ~\textbf{Initial $y$} & ~\textbf{Terminal $y'$} & $\ablkij{k}{R}{yy'}$
      & \textbf{Description}\\\hline
    2 & $\ge 1$ & $(R,I,\io,Z)$ & $(R-1,I,\bt,Z)$  & $R\tht_Z$ & Successful retrial \\\hline
    1 & $\ge 0$ & $(R,I,\io,Z)$ & $(R,I,\io,Z')$  & $q_{ZZ'}$ & Environment (idle) \\
      &         & $(R,I,\io,Z)$ & $(R,I,\bt,Z)$   & $\la_Z$  & Arrival
        while idle\\
      &         & $(R,I,\io,Z)$ & $(R,I,\ga,Z)$     & $\xi_Z$  & Server fails
      while idle \\
      &         & $(R,I,\io,Z)$ & $(R,I,\io,Z)$     & $-\diagi$
        & Diagonal entry (idle) \\
      &         & $(R,I,\bt,Z)$ & $(R,I,\bt,Z')$  & $q_{ZZ'}$ & Environment (busy) \\
      &         & $(R,I,\bt,Z)$ & $(R,I-1,\io,Z)$ & $\mu_Z$ & Demand, $I-1>s$ \\
      &         & $(R,s+1,\bt,Z)$ & $(R,S,\io,Z)$   & $\mu_Z$ & Demand, restocked \\
      &         & $(R,I,\bt,Z)$ & $(R,I,\bt,Z)$   & $-\diagb$
        & Diagonal entry (busy) \\
      &         & $(R,I,\ga,Z)$ & $(R,I,\ga,Z')$      & $q_{ZZ'}$ & Environment (failed) \\
      &         & $(R,I,\ga,Z)$ & $(R,I,\io,Z)$     & $\al_Z$   & Server
        repaired \\
      &         & $(R,I,\ga,Z)$ & $(R,I,\ga,Z)$       & $-\diagg$
        & Diagonal entry (failed) \\\hline
    0 & $\ge 0$ & $(R,I,\bt,Z)$ & $(R+1,I,\ga,Z)$   & $\xi_Z$ & Server fails
      when busy\\
      &         & $(R,I,\bt,Z)$ & $(R+1,I,\bt,Z)$ & $\la_Z$  & Arrival 
        while busy\\
      &         & $(R,I,\ga,Z)$ & $(R+1,I,\ga,Z)$     & $\la_Z$  & Arrival 
        while failed\\\hline
  \end{tabular}\vspace{2ex}
  \caption{Nonzero entries of infinitesimal generator $Q^*$ of
    $\Phi$ \label{tab:A_2^(R)}}  
\end{center}\end{table}

We next formulate $Q^*$ in terms of higher-level block entries. As in Neuts 
\cite{neuts_mgssm_81}, let
$\Dlt(\bv)$ denote the $m\times m$ diagonal matrix whose nonzero entries are the
corresponding entries of the $m$-vector $\bv$. The nonzero
$K$-dimensional square block entries $\ablk{2}{R}$, $\ablk{1}{R}$, and
$\ablk{0}{R}$ of $Q^*$ defined in Eqn. \eqref{eq:Q*} appear as

\begin{align*}\small
  \ablk{2}{R} &= \begin{bmatrix}
    \Tht_R & 0    & 0    & 0    & \cdots & 0 \\
    0    & \Tht_R & 0    & 0    & \cdots & 0 \\
    0    & 0    & \Tht_R & 0    & \cdots & 0 \\
    \vdots & \vdots & \ddots & \vdots & \ddots & \vdots \\
    0    & 0    & 0    & \cdots & \Tht_R & 0 \\
    0    & 0    & 0    & \cdots & 0    & \Tht_R
  \end{bmatrix},\qquad R=1,2,\dots \\[4mm]\small
  \ablk{1}{R} &= \begin{bmatrix}
    \Ga_R  & 0    & 0    & 0    & \cdots & M \\
    M      & \Ga_R  & 0    & 0    & \cdots & 0 \\
    0    & M      & \Ga_R  & 0    & \cdots & 0 \\
    \vdots & \ddots & \ddots & \ddots & \ddots & \vdots \\
    0    & 0    & \dots  & M      & \Ga_R  & 0 \\
    0    & 0    & \dots  & 0    & M      & \Ga_R
  \end{bmatrix}\quad
  \ablk{0}{R} = \begin{bmatrix}
    \La    & 0    & 0    & 0    & \cdots & 0 \\
    0    & \La    & 0    & 0    & \cdots & 0 \\
    0    & 0    & \La    & 0    & \cdots & 0 \\
    \vdots & \ddots & \ddots & \ddots & \ddots & \vdots \\
    0    & 0    & 0    & \cdots & \La    & 0 \\
    0    & 0    & 0    & \cdots & 0    & \La
  \end{bmatrix}, \\[2ex]
  &\hspace{5cm} R=0,1,\dots,
\end{align*}
where the $3m$-dimensional square matrices $\Theta_R$, $\Ga_R$,
$\La$, and $M$ are given by
\begin{align*}
  &\Tht_R =\begin{bmatrix}
    0 & \diagm{R\bth} & 0 \\
    0 & 0           & 0 \\   
    0 & 0           & 0
  \end{bmatrix} \quad
  &
  &\Ga_R = \begin{bmatrix}
    Q_{\io}(R)   & \diagm{\bla} & \diagm{\bxi} \\
    0          & Q_{\bt}      & 0 \\   
    \diagm{\bal} & 0          & Q_{\ga}
  \end{bmatrix} \\[1ex]
  &\La = \begin{bmatrix}
    0 & 0          & 0 \\
    0 & \diagm{\bla} & \diagm{\bxi} \\   
    0 & 0          & \diagm{\bla}
  \end{bmatrix}
  &
  &M = \begin{bmatrix}
    0          & 0 & 0\\
    \diagm{\bmu} & 0 & 0\\
    0          & 0 & 0
  \end{bmatrix}
\end{align*}
and, for each $R\in\Zp$, the $m$-dimensional square matrices $Q_{\io}(R)$,
$Q_{\bt}$, and $Q_{\ga}$ are defined as
\begin{align}\notag
  Q_{\io}(R) &= Q-\diagm{\bla+\bxi+R\bth}, \hspace{2cm} R\in\Zp \\\notag
  Q_{\bt}    &= Q-\diagm{\bla+\bmu+\bxi} \\\label{eq:Q-Delta}
  Q_{\ga}    &= Q-\diagm{\bla+\bal}.
\end{align}
The `0' terms in each of the preceding matrices and those that follow
are square matrices (or scalars) whose dimensions are given by the context
in which they appear.



As a result of the preceding construction, the following result may be stated:

\begin{thm}\label{thm:phi_ldqbd}
  The process $\Phi$ is an irreducible continuous-time LDQBD with infinitesimal
  generator $Q^*$ whose nonzero entries are given in Table \ref{tab:A_2^(R)}.
  Furthermore, the jump process $\wPhi$ of $\Phi$ is a row-convergent LDQBD
  with a transition probability matrix given by $\tilde{P}$ as it appears
  in Eqn. \eqref{eq:P}.
\end{thm}
\begin{proof}
  It remains to show that $\wPhi$ is a row-convergent discrete-time LDQBD. That
  its transition probability matrix is of the form given by Eqn. \eqref{eq:P} is
  a fundamental property of jump processes of continuous-time LDQBDs.
  Accordingly, we begin by constructing the block matrices $\tablk{k}{R}$
  for each $k=0,1,2$ as defined in Eqn. \eqref{eq:P}, followed by the determination of the
  element-wise limit
  \[\tA^*=\lim_{R\to\infty}\tAi{R},\quad\text{where } \tAi{R}=\tablk{0}{R}
    +\tablk{1}{R}+\tablk{2}{R}\]
  if it exists. For convenience, we will define the
  $m$- (row) vectors $\bdiagi$, $\bdiagb$,
  and $\bdiagg$, whose entries consist of terms $\diagi$, $\diagb$, and
  $\diagg$, for each $z=1,\dots, m$. Further, define for each $R\in\Zp$ the
  $(3m)$-dimensional composite block matrix
  \[\Dlt \Ga_R=\begin{bmatrix}
      \Dlt(\bdiagi) & 0 & 0 \\
      0 & \Dlt(\bdiagb) & 0 \\
      0 & 0 & \Dlt(\bdiagg)
    \end{bmatrix}.\]
  We divide each of the rows of $\Ai{R}$
  by the corresponding diagonal (nonzero) entries of $\Dlt \Ga_R$
  to obtain the $(S-s)(3m)$-dimensional square matrix $\tAi{R}$ of the jump
  process: 
  \[\small\tAi{R}= 
    \begin{bmatrix}
      \inv{\Dlt \Ga_R}S_R & 0 & 0 & 0 & \cdots & \inv{\Dlt \Ga_R}M \\
      \inv{\Dlt \Ga_R}M   & \inv{\Dlt \Ga_R}S_R & 0 & 0 & \cdots & 0 \\
      0    & \inv{\Dlt \Ga_R}M & \inv{\Dlt \Ga_R}S_R & 0 & \cdots & 0 \\
      \vdots & \ddots & \ddots & \ddots & \ddots & \vdots \\
      0 & 0 & \dots & \inv{\Dlt \Ga_R}M & \inv{\Dlt \Ga_R}S_R & 0 \\
      0 & 0 & \dots & 0 & \inv{\Dlt \Ga_R}M & \inv{\Dlt \Ga_R}S_R
    \end{bmatrix},\]
  where the $3m$-dimensional square matrix $S_R=\La+\Ga^{\circ}_R+\Tht_R$ and
  \[\Ga^{\circ}_R=\begin{bmatrix}
    Q+\diagm{\bq}   & \diagm{\bla}  & \diagm{\bxi} \\
    0               & Q+\diagm{\bq} & 0 \\   
    \diagm{\bal}    & 0             & Q+\diagm{\bq}
  \end{bmatrix}\]
  is the matrix $\Ga_R$ with scalar diagonal entries set equal to 0.
  
  The subsequent computation of the limiting
  matrix $\tA^*=\lim_{R\to\infty}\tAi{R}$ will be accomplished in an
  element-wise fashion. Its expression will require the $3m\times 3m$
  limiting matrix
  \[\tS^*=\lim_{R\to\infty}{\inv{\Dlt \Ga_R}S_R}=\tLa+\tGa^*+\tTh^*,\]
  for which the terms
  \begin{align}\notag
    \tM    &= \lim_{R\to\infty}\inv{\Dlt \Ga_R}M=\inv{\Dlt \Ga_R}M, &
      \tLa   &= \lim_{R\to\infty}\inv{\Dlt \Ga_R}\La
        =\inv{\Dlt \Ga_R}\La,\\\label{eq:A*_comp}
    \tGa^* &= \lim_{R\to\infty}{\inv{\Dlt \Ga_R}\Ga^{\circ}_R}, &
      \tTh^* &= \lim_{R\to\infty}{\inv{\Dlt \Ga_R}\Tht_R},
  \end{align}
  are evaluated in an element-wise manner. Using the shorthand
  \[\frac{A}{B}=\inv{B}A \quad\text{or}\quad\frac{\trps{A}}{\trps{B}}
    =\trps{A}\trps{\apar{\inv{B}}}\]
  for two square matrices $A$ and $B$, we obtain the limiting
  matrix
  \[\tA^*=\lim_{R\to\infty}\tAi{R}=
    \begin{bmatrix}
      \tS^* & 0 & 0 & 0 & \cdots & \tM \\
      \tM & \tS^* & 0 & 0 & \cdots & 0 \\
      0 & \tM & \tS^* & 0 & \cdots & 0 \\
      \vdots & \ddots & \ddots & \ddots & \ddots & \vdots \\
      0 & 0 & \dots & \tM & \tS^* & 0 \\
      0 & 0 & \dots & 0 & \tM & \tS^*
    \end{bmatrix},\]
  with $3m\times 3m$ block elements given by
  \[\tS^*=
    \begin{bmatrix}
      0 & I_3 & 0  \\
      0 & \frac{\Dlt(\bla+\bq)+Q}{\Dlt(\bdiagb)} & \frac{\Dlt(\bxi) }{\Dlt(\bdiagb)} \\
      \frac{\Dlt(\bal)}{\Dlt(\bdiagg)} & 0 & \frac{\Dlt(\bla+\bq)+Q}{\Dlt(\bdiagg)}
    \end{bmatrix}\qquad
  \text{and}\qquad
  \tM=
    \begin{bmatrix}
      0                              & 0 & 0 \\
      \frac{\Dlt(\bmu)}{\Dlt(\bdiagb)} & 0 & 0 \\
      0                              & 0 & 0
    \end{bmatrix}.\]
  $\wPhi$ is thus a row-convergent LDQBD, which completes the proof of the Theorem.
\end{proof}

\section{System Stability}\label{sec:model_stability}

An analytic traffic intensity formula will now be derived for the inventory model
of this discussion. It is a well-known fact (see Sennott, Humblet, and Tweedie 
\cite{senn83}) that the ergodicity or non-ergodicity of an irreducible
continuous-time LDQBD $\Phi$ is equivalent to that of its embedded, or jump,
chain $\wPhi$. Moreover, as it was shown in Theorem \ref{thm:phi_ldqbd} that $\wPhi$ is
a row-convergent discrete-time LDQBD, the ergodicity condition
Eqn. \eqref{eq:ldqbd_erg} may be used to obtain an analogous drift condition for
its stability, which appears as the following result.
\begin{thm}\label{thm:main}
  The continuous-time LDQBD process $\Phi$ is positive recurrent if and only if
  \begin{equation}\label{eq:pos_rec_crit}
     \bp\cdot\abrk{\bla\bxi+\bal(\bla+\xi)} < \bp\cdot\bal(\bmu+\bxi),
  \end{equation}
  where the $m$-dimensional row vector $\bp$ solves the system of equations given by
  \[\bp Q=\bZ,\quad \bp\eins[m]=1\]
  and $\eins[m]$ is the $m$-dimensional column vector of ones. All multiplicative and
  additive binary relationships in Eqn. \eqref{eq:pos_rec_crit} are performed element-wise,
  save for the operation `$\,\cdot$', which denotes the vector dot product.
\end{thm}\vspace{-3ex}
\begin{proof}
  Let $\wPhi$ be the jump process of $\Phi$. The criterion given in Eqn.
  \eqref{eq:pos_rec_crit} for the positive recurrence of $\Phi$ will be derived
  from the limiting average drift $\tD^*$ of $\wPhi$ that was defined in Eqn. \eqref{eq:D*}.
  In order to compute $\tD^*$, the row-vector solution
  \[\bpist=[\py{i}[j][k]]_{~\shortstack[l]{
      $\sss i = s+1,\dots,S$ \\
      $\sss j = \io,\bt,\ga$ \\
      $\sss k = 1,\dots,m$}}\]
  of the system expressed by Eqn. \eqref{eq:pi*_system} is required. Note that the vector
  is written in partitioned form according to the states $(i,j,k)\in\Sphi$. For example, the
  notation $\py{i}$ denotes the $3m$-dimensional partition of
  $\bpist$ for which $i$ is held constant and $\py{i}[j]$ the
  $m$-dimensional partition for which both $i$ and
  $j$ are fixed.
  
  When expanded, the system of equations expressed by Eqn. \eqref{eq:pi*_system} becomes
   \begin{align}\small\notag
     \py{(s+1)}\tS^*  + \py{(s+2)}\tM  &=  \py{(s+1)} \\\notag
     \py{(s+2)}\tS^*  + \py{(s+3)}\tM  &=  \py{(s+2)} \\\notag
       &\mathrel{\makebox[\widthof{=}]{\vdots}} \\\notag
     \py{(s+1)}\tM    +  \py{S}\tS^*   &=  \py{S} \\\label{eq:expanded_pi_sys}
     \sum_{i=s+1}^S\py{i}\eins[3m]     &=  1.
  \end{align}

  We will proceed by induction on the inventory difference term $(S-s)$.
  Consider an inventory system in which $S-s=2$.
  The linear system in Eqn. \eqref{eq:expanded_pi_sys} may be written in
  vector-matrix form as
  \begin{equation}\small\label{eq:eqn9rewritten}
    \begin{bmatrix}
      \py{(s+1)} & \py{(s+2)}
    \end{bmatrix}
    \begin{bmatrix}
      \tS^* & \tM   & \eins[3m] \\
      \tM   & \tS^* & \eins[3m] \\
    \end{bmatrix} =
    \begin{bmatrix}
      \py{(s+1)} & \py{(s+2)} & 1
    \end{bmatrix}
  \end{equation}  
with the partitioned vector solution.
  \[\bpist_2=
    \begin{bmatrix}
      \py{(s+1)} & \py{(s+2)}
    \end{bmatrix}.\]
  For convenience, we will now write Eqn. \eqref{eq:eqn9rewritten}
  as the transpose system
  \[\small\begin{bmatrix}
      \trps{(\tS^*)}   & \trps{(\tM)}    \\
      \trps{(\tM)}     & \trps{(\tS^*)}  \\
      \trps{(\eins[3m])} & \trps{(\eins[3m])}
    \end{bmatrix}
    \begin{bmatrix}
      \trps{\py{(s+1)}} \\ \trps{\py{(s+2)}}
    \end{bmatrix} =
    \begin{bmatrix}
      \trps{\py{(s+1)}} \\ \trps{\py{(s+2)}} \\ 1
    \end{bmatrix}\]
which, when expanded, becomes
  \begin{align}\notag &\small
    \begin{bmatrix}
      0 & 0 &  \frac{\Dlt(\bal)}{\Dlt(\bdiagg)}
        & 0 & \frac{\Dlt(\bmu) }{\Dlt(\bdiagb)} & 0 \\
      I_3 & \frac{\Dlt(\bla+\bq)+\trps{Q}}{\Dlt(\bdiagb)} & 0 & 0 & 0 & 0 \\
      0 & \frac{\Dlt(\bxi)}{\Dlt(\bdiagb)} & \frac{\Dlt(\bla+\bq)+\trps{Q}}{\Dlt(\bdiagg)}
        & 0 & 0 & 0 \\
      0 & \frac{\Dlt(\bmu) }{\Dlt(\bdiagb)} & 0 
        & 0 & 0 &  \frac{\Dlt(\bal)}{\Dlt(\bdiagg)} \\
      0 & 0 & 0 & I_3 & \frac{\Dlt(\bla+\bq)+\trps{Q}}{\Dlt(\bdiagb)} & 0 \\
      0 & 0 & 0 & 0
        & \frac{\Dlt(\bxi)}{\Dlt(\bdiagb)} & \frac{\Dlt(\bla+\bq)+\trps{Q}}{\Dlt(\bdiagg)} \\
      \trps{(\eins[m])} & \trps{(\eins[m])} & \trps{(\eins[m])} & \trps{(\eins[m])}
        & \trps{(\eins[m])} & \trps{(\eins[m])} \\
    \end{bmatrix}
    \begin{bmatrix}
      \trps{\py{(s+1)}[\io]} \\ \trps{\py{(s+1)}[\bt]} \\ \trps{\py{(s+1)}[\ga]} \\
        \trps{\py{(s+2)}[\io]} \\ \trps{\py{(s+2)}[\bt]} \\ \trps{\py{(s+2)}[\ga]}
    \end{bmatrix} \\[3ex]\small \label{eq:exp_blk_sys} &=
    \trps{\begin{bmatrix}
      \py{(s+1)}[\io] & \py{(s+1)}[\bt] & \py{(s+1)}[\ga] &
      \py{(s+2)}[\io] & \py{(s+2)}[\bt] & \py{(s+2)}[\ga] & 1
    \end{bmatrix}},
  \end{align}
  where each occurrence of the symbol 0 represents the $m\times m$ array of zeroes.
  The solution $\bpist_2$ of Eqn. \eqref{eq:exp_blk_sys}, re-expressed as a column vector with
  $m$-entry partitions, is most expediently obtained if one first solves the system in terms of
  the block matrix entries of the corresponding coefficient matrix, from whence we will obtain
  a vector solution with block matrix entries. This solution may then be easily converted to
  the requisite scalar vector solution of Eqn. \eqref{eq:exp_blk_sys}.
  
  Variable substitutions will now be made in Eqn. \eqref{eq:exp_blk_sys} in
  order to accommodate the expression of systems of $m\times m$ matrix terms. First,
  we replace each occurrence of the $m$-dimensional vector variables
  $\py{I}[x]$ with the $m\times m$ matrix terms $\pyblk{I}{x}$, for each $I\in\Zp$ and for
  each $x=\io,\bt,\ga$. Consequently, one may
  perform the conversion from the vector with diagonal matrix entries to the corresponding
  $m$-entry row vector via the relationship
  \begin{equation}\label{eq:v_block_to_v_id}
    \py{I}[x]=\trps{(\eins[m])}\cdot\pyblk{I}{x}.
  \end{equation}

  Next, we replace the entries of the last row of the coefficient matrix with
  the identity matrix $I_m$. For the moment, the $1$ on the right-hand side
  will be replaced with the indeterminate quantity $\bu\in\Rp^{m\times m}$ until an
  appropriate value can be determined.

  We may then rewrite the system in Eqn. \eqref{eq:exp_blk_sys} as
  \begin{align}\notag &\small
    \begin{bmatrix}
      0 & 0 & \frac{\Dlt(\bal)}{\Dlt(\bdiagg)}
        & 0 & \frac{\Dlt(\bmu) }{\Dlt(\bdiagb)} & 0 \\
      I_m & \frac{\Dlt(\bla+\bq)+\trps{Q}}{\Dlt(\bdiagb)} & 0 & 0 & 0 & 0 \\
        0 & \frac{\Dlt(\bxi)}{\Dlt(\bdiagb)} & \frac{\Dlt(\bla+\bq)+\trps{Q}}{\Dlt(\bdiagg)}
          & 0 & 0 & 0 \\
      0 & \frac{\Dlt(\bmu) }{\Dlt(\bdiagb)} & 0 
        & 0 & 0 &  \frac{\Dlt(\bal)}{\Dlt(\bdiagg)} \\
      0 & 0 & 0 & I_m & \frac{\Dlt(\bla+\bq)+\trps{Q}}{\Dlt(\bdiagb)} & 0 \\
      0 & 0 & 0 & 0
        & \frac{\Dlt(\bxi)}{\Dlt(\bdiagb)} & \frac{\Dlt(\bla+\bq)+\trps{Q}}{\Dlt(\bdiagg)} \\
      I_m & I_m & I_m & I_m & I_m & I_m \\
    \end{bmatrix}
    \begin{bmatrix}
      \trps{\pyblk{(s+1)}{\io}} \\ \trps{\pyblk{(s+1)}{\bt}} \\ \trps{\pyblk{(s+1)}{\ga}} \\
      \trps{\pyblk{(s+2)}{\io}} \\ \trps{\pyblk{(s+2)}{\bt}} \\ \trps{\pyblk{(s+2)}{\ga}}
    \end{bmatrix} \\[3ex] \label{eq:exp_blk_sys2} &=
    \trps{\begin{bmatrix}
      \pyblk{(s+1)}{\io} & \pyblk{(s+1)}{\bt} & \pyblk{(s+1)}{\ga} &
      \pyblk{(s+2)}{\io} & \pyblk{(s+2)}{\bt} & \pyblk{(s+2)}{\ga} & \bu
    \end{bmatrix}}.
  \end{align}
  We denote the solution of the system in Eqn. \eqref{eq:exp_blk_sys2} as 
  the (row) vector $\bpist_b$ of $m\times m$ matrix entries. Once it has been established that
  the individual entries $\bpist_{I\eta}$ of $\bpist_b$,
  for $I=s+1,\,s+2$ and $\eta\in\sett{\io,\bt,\ga}$,
  of this system are diagonal matrices, then we may, in a manner analogous to that of
  Eqn. \eqref{eq:v_block_to_v_id}, say that the vector-multiplicative operation
  \begin{equation}\label{eq:v_block_to_v_id_2}
    \begin{bmatrix}
      \trps{(\eins[m])} & \trps{(\eins[m])} & \trps{(\eins[m])} & \trps{(\eins[m])} &
      \trps{(\eins[m])} & \trps{(\eins[m])}
    \end{bmatrix}\cdot\bpist_b
    =\trps{(\eins[K])}\cdot\bpist_b
  \end{equation}
  yields a $K$-dimensional row vector with scalar entries.

  In order to solve Eqn. \eqref{eq:exp_blk_sys2} using conventional methods for
  linear systems with scalar unknowns, it would be necessary for all elements of the coefficient matrix to be diagonal matrices. However, the generator matrix $Q$
  of the random environment is not diagonal. In addition, since row sums of $A^*$ are not multiples of $I_m$, it is not `stochastic' in the block-matrix sense, which, in effect, causes the system to become inconsistent for any value of $\bu$. To overcome this difficulty, we will first solve for what will be termed a $Q$-homogeneous solution $\bpist_h$ of Eqn. \eqref{eq:exp_blk_sys2} by setting $Q=0$, which then results in the matrix $A^*$ becoming `stochastic' in the block-matrix sense. We thus solve

  \begin{equation}\label{eq:exp_blk_sys3}\small
    \begin{bmatrix}
      0 & 0 &  \frac{\Dlt(\bal)}{\Dlt(\bdiagg)}
        & 0 & \frac{\Dlt(\bmu) }{\Dlt(\bdiagb)} & 0 \\
      I_m & \frac{\Dlt(\bla+\bq)}{\Dlt(\bdiagb)} & 0 & 0 & 0 & 0 \\
      0 & \frac{\Dlt(\bxi)}{\Dlt(\bdiagb)} & \frac{\Dlt(\bla+\bq)}{\Dlt(\bdiagg)}
        & 0 & 0 & 0 \\
      0 & \frac{\Dlt(\bmu) }{\Dlt(\bdiagb)} & 0 
        & 0 & 0 &  \frac{\Dlt(\bal)}{\Dlt(\bdiagg)} \\
      0 & 0 & 0 & I_m & \frac{\Dlt(\bla+\bq)}{\Dlt(\bdiagb)} & 0 \\
      0 & 0 & 0 & 0
        & \frac{\Dlt(\bxi)}{\Dlt(\bdiagb)} & \frac{\Dlt(\bla+\bq)}{\Dlt(\bdiagg)} \\
      I_m & I_m & I_m & I_m & I_m & I_m \\
    \end{bmatrix}
    \begin{bmatrix}
      \trps{\pyblkh{(s+1)}{\io}} \\ \trps{\pyblkh{(s+1)}{\bt}} \\
      \trps{\pyblkh{(s+1)}{\ga}} \\ \trps{\pyblkh{(s+2)}{\io}} \\
      \trps{\pyblkh{(s+2)}{\bt}} \\ \trps{\pyblkh{(s+2)}{\ga}}
    \end{bmatrix}
    =
    \begin{bmatrix}
      \trps{\pyblkh{(s+1)}{\io}} \\ \trps{\pyblkh{(s+1)}{\bt}} \\
      \trps{\pyblkh{(s+1)}{\ga}} \\ \trps{\pyblkh{(s+2)}{\io}} \\
      \trps{\pyblkh{(s+2)}{\bt}} \\ \trps{\pyblkh{(s+2)}{\ga}} \\ \bu_h
    \end{bmatrix},
  \end{equation}
  where we allow $\bu_h\in\Rp^{m\times m}$ to be an indeterminate quantity. Using a
  symbolic linear equation solver, we thus obtain
  \[\small\trps{(\bpist_h)}=\frac{1}{2}\cdot
    \small
    \begin{bmatrix}
      \Dlt(\bal)\Dlt(\bmu+\bxi) \\
      \Dlt(\bal)\Dlt(\bq+\bla+\bmu+\bxi) \\
      \Dlt(\bxi)\Dlt(\bq+\bla+\bal) \\
      \Dlt(\bal)\Dlt(\bmu+\bxi) \\
      \Dlt(\bal)\Dlt(\bq+\bla+\bmu+\bxi) \\
      \Dlt(\bxi)\Dlt(\bq+\bla+\bal)
    \end{bmatrix},\]
  where
  \[\bu_h=\Dlt(\bal)\Dlt(2\bmu+2\bxi+\bla+\bq)+\Dlt(\bxi)\Dlt(\bq+\bla+\bal).\]

  Next, the non-$Q$-homogeneous system in Eqn. \eqref{eq:exp_blk_sys2} will be solved.
  In order to do this, we first define the $m\times m$ matrix
  \[\Pi=\trps{
    \begin{bmatrix}
      \trps{\bp} & \trps{\bp} & \dots & \trps{\bp}
    \end{bmatrix}},\]
  where $\bp$ is the stationary probability vector of the random
  environment that was defined in the statement of the Theorem.

  \begin{prop} \label{prop:sys12sol}
    The unique block matrix solution to the linear system in
    Eqn. \eqref{eq:exp_blk_sys2} with $\bu=\bu_h\cdot\trps{\Pi}$ is
    \[\small\trps{(\bpist_b)}=\trps{(\bpist_h)}\cdot\trps{\Pi}.\]
    This establishes that the entries
    $\tbpi_{I\eta}$ of $\bpist_b$, for $I=s+1,\,s+2$ and
    $\eta\in\sett{\io,\bt,\ga}$, are diagonal square matrices,
    and thus that the relations given in Eqn. \eqref{eq:v_block_to_v_id} and
    Eqn. \eqref{eq:v_block_to_v_id_2} are valid.
  \end{prop}
  
  \begin{proof} (Proposition \ref{prop:sys12sol})
    That $\trps{(\bpist_b)}$ is a solution to Eqn. \eqref{eq:exp_blk_sys2} may be
    verified by evaluating the system in Eqn. \eqref{eq:exp_blk_sys2} with the given
    value of $\trps{(\bpist_b)}$ and subsequently applying the identity
    \[\bp\cdot Q=\trps{Q}\trps{\bp}=\bZ.\]
    Uniqueness is a consequence of the fact that $\trps{(\bpist_b)}$ is a
   `stationary vector' of the system in Eqn. \eqref{eq:exp_blk_sys2}.
  
    Finally, we will state without formal demonstration that, regardless
    of whether one evaluates the system in Eqn. \eqref{eq:exp_blk_sys2} as a
    block-matrix system with the column vector solution $\trps{(\bpist_b)}$ of
    $m$-dimensional matrix entries or as a scalar system with the $K=3m(S-s)$
    scalar-entry row vector solution $\trps{(\eins[K])}\bpist_b$,
    equivalent results are produced (up to a block-matrix interpretation).
    This is a consequence of the fact that $\bpist_b$ is expressed entirely in
    terms of diagonal matrices.
  \end{proof}

  We may use Proposition 1 to construct a vector solution
  $\bpist_2$ with scalar entries to the system in Eqn. \eqref{eq:exp_blk_sys},
  which is detailed in the following result:

  \begin{prop} \label{prop:sys10sol}
    The $K$-column vector solution
    $\trps{(\bpist_2)}$ to the linear system in Eqn. \eqref{eq:exp_blk_sys}
    with scalar entries is
    \[\trps{(\bpist_2)}=\frac{1}{mc}\trps{(\bpist_b)}\cdot\eins[K]\]
    where
    \begin{equation}\label{eq:c}
      c=\bp\cdot\abrk{\bal(2\bmu+2\bxi+\bla+\bq)+\bxi(\bq+\bla+\bal)}.
    \end{equation}
    The operation `$\cdot$' in Eqn \eqref{eq:c} is defined as the vector dot
    product while other binary operations between vectors are performed element-wise.
  \end{prop}

  \begin{proof} (Proposition \ref{prop:sys10sol}) 
    It is first necessary to apply Eqn. \eqref{eq:v_block_to_v_id} in order to
    convert $\bpist_b$ into a (column) $K$-vector term, which is then
    normalized into a probability vector through division by the following scalar:
    \begin{align*}
      \trps{(\eins[K])}\bpist_b\eins[K] &= \trps{(\eins[m])}\trps{\bu}\eins[m] \\
                        &=    \trps{(\eins[m])}\apar{\Pi\cdot\bu_h}\eins[m] \\
                        &=     mc.
    \end{align*}
    Substituting the
    resulting expression, defined as $\bpist_2$ in the statement of Proposition
    \ref{prop:sys10sol}, into the linear system in Eqn. \eqref{eq:exp_blk_sys} shows
    that $\bpist_2$ is indeed a solution to this system.
\end{proof}

Now that a limiting stationary vector $\bpist_2$ is in hand, we proceed to
compute the corresponding limiting drift expression. First,
we observe that $\bpist_2$ is composed of repeating blocks of 
$3m$-dimensional vectors $\trps{(\bpi_r)}$, where
\[\small\trps{(\bpi_r)}=
    \begin{bmatrix}
      \Dlt(\bal)\Dlt(\bmu+\bxi) \\
      \Dlt(\bal)\Dlt(\bq+\bla+\bmu+\bxi) \\
      \Dlt(\bxi)\Dlt(\bq+\bla+\bal) \\
    \end{bmatrix},\]
which yields 
\[\bpist_2=\frac{1}{2mc}\abrk{\trps{(\eins[m])}\Pi\cdot\bpi_r\quad
  \trps{(\eins[m])}\Pi\cdot\bpi_r}.\]
  By Eqn. \eqref{eq:D*}, we compute
  \begin{align}\notag\small
    \tD^* &= \bpist_2(\tA_0^* - \tA_2^*)\eins \\\notag
        &= \frac{1}{2mc}\trps{(\eins[m])}\Pi\cdot\abrk{\bpi_r\quad\bpi_r}\cdot
           \begin{bmatrix}
             \tLa - \tTh^* & 0 \\
             0 & \tLa - \tTh^*
           \end{bmatrix}
           \begin{bmatrix}
            \eins[3m] \\ \eins[3m]
           \end{bmatrix}\\\notag
        &= \frac{1}{mc}\trps{(\eins[m])}\Pi\cdot\bpi_r\cdot(\tLa - \tTh^*)
           \eins[m] \\\notag
        &= \frac{1}{mc}\trps{(\eins[m])}\Pi\cdot\bpi_r\cdot
           \begin{bmatrix}
              0 & -I_m & 0  \\
              0 & \frac{\Dlt(\bla)}{\Dlt(\bq+\bla+\bmu+\bxi)}
                & \frac{\Dlt(\bxi)}{\Dlt(\bq+\bla+\bmu+\bxi)} \\
              0 & 0 & \frac{\Dlt(\bla)}{\Dlt(\bq+\bla+\bal)}
           \end{bmatrix}
           \begin{bmatrix}
            \eins[m] \\ \eins[m] \\ \eins[m]
           \end{bmatrix} \\\notag
        &= \frac{1}{mc}\trps{(\eins[m])}\cdot
           \begin{bmatrix}
             0 & \Pi\Dlt(\bla)\Dlt(\bal)-\Pi\Dlt(\bal)\Dlt(\bmu+\bxi)
               & \Pi\Dlt(\bxi)\Dlt(\bla+\bal)
           \end{bmatrix}           
           \begin{bmatrix}
            \eins[m] \\ \eins[m] \\ \eins[m]
           \end{bmatrix} \\\notag
        &= \frac{1}{mc}
            \abrk{m\bp\cdot\abrk{\bla\bal-\bal(\bmu+\bxi)+\bxi(\bla+\bal)}}
            \\[2ex]
           \label{eq:(S-s)=2_drift}
        &= \frac{1}{c}\bp\cdot
           \abrk{\bla\bxi+\bal(\bla+\xi)-\bal(\bmu+\bxi)}.
  \end{align}

  For the induction step, we assume that the drift expression
  Eqn. \eqref{eq:(S-s)=2_drift} holds for $(S-s-1)$. The stationary probability
  vector for this model may then be obtained as
  \[\bpist_{(S-s-1)}=\frac{1}{(S-s-1)mc}\underbrace{\abrk{\trps{(\eins[m])}\Pi\cdot\bpi_r\quad
    \trps{(\eins[m])}\Pi\cdot\bpi_r\quad\dots\quad\trps{(\eins[m])}\Pi\cdot\bpi_r}}_{(S-s-1)
    \text{ terms}}.\]
  For an $(S-s)$ model, the matrix $\tA^*$ gains an additional repeated block
  matrix row, from which we deduce the new stationary probability vector
  to be
  \[\bpist_{(S-s)}=\frac{1}{(S-s)mc}\underbrace{\abrk{\trps{(\eins[m])}\Pi\cdot\bpi_r\quad
    \trps{(\eins[m])}\Pi\cdot\bpi_r\quad\dots\quad\trps{(\eins[m])}\Pi\cdot\bpi_r}}_{(S-s)
    \text{ terms}}.\]
  We now repeat the previous computation of drift $D^*$ as  
  \begin{align}\notag
    \tD^* &= \bpist_{(S-s)}(\tA_0^* - \tA_2^*)\eins \\\notag
        &= \frac{1}{(S-s)mc}\trps{(\eins[m])}\Pi\abrk{\bpi_r\quad\bpi_r
             \quad\dots\quad\bpi_r} \\[2ex]\notag
        &\small\hspace{2em}\times
          \begin{bmatrix}
            \tLa - \tTh^* & 0 & 0 & 0 & \dots & 0 \\
            0 & \tLa - \tTh^* & 0 & 0 & \dots & 0 \\
            0 & 0 & \tLa - \tTh^* & 0 & \dots & 0 \\
            \vdots & \vdots & \ddots & \ddots &  & \vdots \\
            0 & 0 & 0 & 0 & \dots & \tLa - \tTh^*
          \end{bmatrix}
          \begin{bmatrix}
           \eins[m] \\ \eins[m] \\ \eins[m] \\ \vdots \\ \eins[m]
          \end{bmatrix} \\[2ex]\notag
        &= \frac{1}{mc}\trps{(\eins[m])}\Pi\bpi_r(\tLa - \tTh^*)\eins[m]
          \\[2ex]\notag
        &= \frac{1}{mc}\trps{(\eins[m])}\Pi\bpi_r
           \begin{bmatrix}
              0 & -I_m & 0  \\
              0 & \frac{\Dlt(\bla)}{\Dlt(\bq+\bla+\bmu+\bxi)}
                & \frac{\Dlt(\bxi)}{\Dlt(\bq+\bla+\bmu+\bxi)} \\
              0 & 0 & \frac{\Dlt(\bla)}{\Dlt(\bq+\bla+\bal)}
           \end{bmatrix}
           \begin{bmatrix}
            \eins[m] \\ \eins[m] \\ \eins[m]
           \end{bmatrix} \\[2ex] \label{eq:D*_final}
        &= \frac{1}{c}\bp\cdot
           \abrk{\bla\bxi+\bal(\bla+\xi)-\bal(\bmu+\bxi)}.
  \end{align}
  By then setting $\tD^*<0$, we obtain the expression in Eqn. \eqref{eq:pos_rec_crit}
  for the positive recurrence of $\Phi$, and the Theorem is proven.
\end{proof}


By reformulating the average drift $\tD^*$ in Eqn. \eqref{eq:D*_final} as
a traffic intensity, the performance measure of average server occupancy of $\Phi$
in steady state is obtained. This is accomplished by setting $\tD^*<0$ and rearranging terms, 
which leads to the following Corollary to Theorem \ref{thm:main}.
\begin{cor}\label{cor:traff_int}
  The traffic intensity $\rho$ of the process $\Phi$ may be written as
  \begin{equation}\label{eq:rho}
    \rho=\frac{\bp\cdot\abrk{\bla\bxi+\bal(\bla+\xi)}}
    {\bp\cdot\bal(\bmu+\bxi)}.
  \end{equation}
  Subsequently, the continuous-time LDQBD process $\Phi$ is
  positive recurrent if and only if $\rho<1$.
\end{cor}

\section{Steady-State Distribution and Performance Measures}\label{sec:measures}

If $D^*<0$, then by Theorem \ref{cor:traff_int}, $\Phi$ is positive
recurrent. In this case, the joint steady-state probabilities are defined as
\begin{align*}
  \prb{R,i,x,z}&=\lim_{t\to\infty}\prb{(R(t),I(t),X(t),Z(t))=(R,i,x,z)}, \\
    &\hspace{3cm} R\in\Zp,\;(i,x,z)\in\mcL
\end{align*}
exist. Since $\mcL$ is a finite set, we may enumerate the elements of this set
as $\mcL=\sett{1,2,\dots,M}$, where we define the $k$th element of $\mcL$ as
$(i_k,x_k,z_k)$ and $M=(S-s)\cdot(3m)$. The steady-state probabilities
may then be expressed more concisely as
\[\prR{R,k}=\prb{R,i_k,x_k,z_k},\]
whereupon we may define the $M$-dimensional row vectors
\[\bprR{R} = (\prR{R,1},\prR{R,2},\dots,\prR{R,M}),\qquad R\in\Zp\]
of steady-state probabilities of $\Phi$ grouped according to orbit size
$R$. Assuming the positive recurrence of $\Phi$,
one may infer the presence of the matrix-geometric relationship between
terms of $\bprR{R}$, which is given in Bright and Taylor \cite{bri_tay95} for
$R\ge 1$ as
\begin{equation}\label{eq:bp_R}
  \bprR{R}=\bprR{0}\sum_{R=1}^{\infty}\abrk{\prod_{\ell=0}^{R-1} R_{\ell}}e,
\end{equation}
where the rate matrices $\sett{R_{\ell}}[\ell\in\Zp]$ are the minimal
non-negative solutions to the system of equations
\begin{equation}\label{eq:R_k}
  \ablk{\ell}{0}+R_{\ell}\ablk{\ell+1}{1}+R_{\ell}
  \abrk{R_{\ell+1}\ablk{\ell+2}{1}}=0,\qquad \ell\in\Zp,
\end{equation}
and the level 0 steady state probability $\bp_0$ is the minimal vector
solution to
\begin{equation}\label{eq:p_0}
  \bprR{0}\abrk{\ablk{1}{0}+R_0\ablk{2}{1}}=0.
\end{equation}

Since it is unlikely that Eqn. \eqref{eq:R_k} and Eqn. \eqref{eq:p_0} have closed-form solutions; however, it is more expedient to produce estimated measures of performance. To this end, one may apply one of several established algorithms that were developed for the purpose of estimating the steady-state distribution of an LDQBD,
such as that of Bright and Taylor \cite{bri_tay95}. The method, via Algorithm 1, produces estimated stationary probabilities $p_{_{R,k}}(R^*)\sim p_{_{R,k}}$
of a truncated system $\Phi(R^*)$, say, at some level (orbit size)
$R^*$ that is sufficiently large. The term `sufficiently large' is
used in the context of the fact that
\[\bprR{R}=\lim_{R^*\to\infty}\bprR{R}(R^*),\quad k\in\Zp.\]
In other words, the estimates become progressively more accurate as the system is truncated at larger levels $R^*$. Because Algorithm 1 produces successive estimates of $\bp_R$ by means of the matrix geometric recurrence relation Eqn. \eqref{eq:bp_R},
there is a need to efficiently compute the rate matrices $R_{\ell}$, a task for which Algorithm 2 is utilized.

With the steady-state distribution of the system $\Phi$ in hand,
the asymptotic performance measures of the queueing inventory system
$\Phi$ may be obtained, beginning with the marginal steady-state
probabilities $p_i,\,p_b,\,p_f$ of the server status:
\begin{align*}
   &\text{Idle Probability:} &
    \prR{\io}&=\sum_{R=0}^{\infty}\sum_{i=s+1}^{S}\sum_{k=1}^m \prb{R,i,\io,k} \\
   &\text{Busy Probability:} &
    \prR{\bt}&=\sum_{R=0}^{\infty}\sum_{i=s+1}^{S}\sum_{k=1}^m \prb{R,i,\bt,k} \\
   &\text{Failure Probability:} &
    \prR{\ga}&=\sum_{R=0}^{\infty}\sum_{i=s+1}^{S}\sum_{k=1}^m \prb{R,i,\ga,k} 
\end{align*}
Likewise, the steady-state probability $\prR{R}$ of the number of
demands in orbit is the marginal probability
\[\prR{R} = \sum_{i=s+1}^{S}\sum_{k=1}^m \abrk{\prb{R,i,\io,k}
  +\prb{R,i,\bt,k}+\prb{R,i,\ga,k}}\]
The long-run expected number of demands  in orbit ($L_R$) and
the system (L) may then be expressed in the usual way as
\[L_R=\sum_{R=0}^{\infty} R\cdot \prR{R} \quad\text{and}\quad L=L_R + \prR{\bt}.\]

Temporal measures of queueing performance require the long-run average
exponential input rate over environment states, which for the stationary
probability vector $\bp=[p_z]$ of $Q$, is given by
\[\barla=\sum_{z=1}^m\la_z\cdot p_z=\bla\cdot\bp.\]
We may then apply Little's Law to obtain the long-run expected wait times of demands in orbit
($W_R)$ and in the system ($W$):
\begin{equation}\label{eq:little_law}
  W_R=L_R/\barla\quad\text{and}\quad W=L/\barla.
\end{equation}

The next result provides for the independence of all performance measures defined thus far on the state $I$ of the inventory.
\begin{prop}\label{prop:ind_of_perf_meas}
  The performance measures $L$, $L_R$, $W$, and $W_R$ are independent of the
  inventory threshold values $s$ and $S$.
\end{prop}
\begin{proof}
  As may be ascertained from the nonzero blocks $\ablk{j}{i}$, $j=0,1,2$, of the
  infinitesimal generator matrix $Q^*$ whose entries are listed in Table
  \ref{tab:A_2^(R)}, the evolution of the inventory state $I$ does not affect
  any of the exponential rates that appear in the third column of the Table, either
  through the appearance of $I$ or of any rate term that pertains to restocking
  delay or some other duration related to the number in inventory.
  Thus, the marginal distributions $\sett{p_R}[R\in\Zp]$ of orbit size, as well
  as the probabilities $\sett{p_{\io}, p_{\bt}, p_{\xi}}$ of server state, along
  with any performance measures derived from these probabilities, do not
  depend on the value of $I$, and hence of $s$ or $S$.
\end{proof}

The system performance-measures of inventory level, depletion, and replenishment
require the steady-state distribution of the number of products in the inventory,
which is given by the marginal long-run probabilities of there being $I\in[s+1,S]$ in inventory
\[\prR{I} = \sum_{R=0}^{\infty}\sum_{k=1}^m \abrk{\prb{R,i,\io,k}
  +\prb{R,i,\bt,k}+\prb{R,i,\ga,k}}.\]
From this distribution, we may obtain the long-run expected inventory level as
\begin{equation}\label{eq:Binv}
  \binv=\sum_{I=s+1}^S{I\cdot\prR{I}}.
\end{equation}
The computation of $\binv$ is greatly simplified by the fact
that its value is dependent solely upon the quantities $s$ and $S$,
as stated and proved in the following Lemma.
\begin{lemma}\label{lemma:B_inv}
  The steady-state probability distribution of the amount in inventory
  for the positive recurrent process $\Phi$ is given by
  \[\prR{I}=\frac{1}{S-s},\quad I=s+1,\dots,S.\]
  Consequently, the expected inventory content $\binv$ may be computed as
  \begin{equation}\label{eq:newBinv}
    \binv=\frac{1}{2}(s+S+1).
  \end{equation}
\end{lemma}
\begin{proof}
  See the Appendix.
\end{proof}

For the long-run expected time $D_S$ to deplete (or replenish) the inventory from the maximum level $S$, we observe that the inventory decrements by one just before a demand
exits the system. Thus, depletion from the maximum level of $S$ items occurs whenever $S-s$ customers are successfully processed, which is, on average, $(S-s)$ average system sojourn times $W$. Thus,
\[D_S = (S-s)W.\]
To obtain the long-run ordering rate $r_o$, we use the fact that
there is one order per depletion time $D_s$ so that
\[r_o=\frac{1}{D_S}=\frac{1}{(S-s)W}.\]
On the other hand, the long-run supply rate $r_s$ is given by the
number of items ordered per depletion time $D_S$. Thus,
\[r_s=\frac{(S-s)}{D_S}=(S-s)r_o=\frac{1}{W}.\]

\section{Optimization Study}\label{sec:optimization}

In this section, the minimal long-run average costs for the operation of three stable inventory
systems of the type described in Section \ref{sec:model} are considered. The objective here is to
compare and contrast the optimal inventory threshold parameters $s$ and $S$ that correspond to the
minimum operational cost of systems in steady state over increasing traffic intensity $\rho$. In what
follows, the construction of the inventory systems of interest, along with the formulation of the
steady state cost function from queueing parameters and associated steady state performance measures,
is described. A method is then given to determine a unique minimal-cost pair $\sSstar$
(up to a choice of $s^*\in\Zp$) for any inventory system of the type described in this paper.

\subsection{System Definitions\label{ss:sys_defn}}

The construction of stable inventory systems of varying traffic intensity may be accomplished
through appropriate choices of the exponential parameters $\la_z$, $\mu_z$, $\xi_z$, and $\al_z$ that
produce increasing values of $\rho$ within the interval $[0,1]$. In addition,
it is ensured that several states $z\in\mcS$ in each system
exhibit values of the single-environment traffic intensity function $\rho_z$ (derived
from Eqn. \eqref{eq:rho} with $\bp=[1]\equiv 1$), given by
\[\small \rho_z = \frac{\la_z+\xi_z}{\mu_z+\xi_z}
   + \frac{\la_z\xi_z}{\al_z(\mu_z+\xi_z)},\]
that are greater than 1, despite having an overall
traffic intensity of $\rho<1$. 
The resulting exponential parameters for each system, together with
$\rho_z$ for each $z\in\mcS$, appear in Tables
\ref{tab:params_lt}, \ref{tab:params_mt}, and \ref{tab:params_ht},
respectively.
\begin{table}[ht]\centering
    \caption{Parameter values of the low-traffic system ($\rho=0.17$) with seven
      environment states.
      \label{tab:params_lt}}
    \begin{tabular}{c|r*{7}{c}|}\hline
      Environment ($z$) & $\la_z$ & $\mu_z$ & $\xi_z$ & $\al_z$
      & $\tht_z$ & $\rho_z$ \\\hline
      1 &   1.0  & 13.0  & 0.05   & 7.0  & 1.00 & 0.0810 \\
      2 &   8.0 &  1.2  & 3.80   & 0.8  & 0.10 & 9.9600\\
      3 &   0.3  & 17.0  & 0.02   & 15.0 & 4.00 & 0.0188\\
      4 &   2.0  & 12.0  & 0.30   & 12.0 & 2.00 & 0.1911\\
      5 &   0.5  & 18.7  & 1.00   & 5.0  & 5.00 & 0.0812\\
      6 &   1.0  & 15.0  & 1.20   & 2.8  & 0.10 & 0.1623 \\
      7 &   5.0 & 6.0   &  4.00  & 0.5  & 0.05 & 4.9000\\\hline
    \end{tabular}
\end{table}
\begin{table}[ht]\centering
    \caption{Parameter values of the medium-traffic system ($\rho=0.41$) with seven
      environment states.
      \label{tab:params_mt}}
    \begin{tabular}{c|r*{7}{c}|}\hline
      Environment ($z$) & $\la_z$ & $\mu_z$ & $\xi_z$ & $\al_z$
      & $\tht_z$ & $\rho_z$\\\hline
      1 &  6.0  & 7.0 & 5.0 & 2.0  &  1.0 & 2.1667\\
      2 &  0.1  & 4.2 & 0.8 & 0.1  &  1.7 & 0.3400\\
      3 &  1.0  & 8.0 & 1.0 & 15.0 &  2.0 & 0.2296\\
      4 &  0.8  & 10.0& 0.3 & 12.0 &  2.0 & 0.1087\\
      5 &  2.0  & 4.5 & 1.0 &  5.0 &  5.0 & 0.6182\\
      6 &  0.5  & 2.0 & 0.7 & 13.0 &  1.0 & 0.4544\\
      7 &  9.0  & 0.3 & 0.2 &  0.5 &  0.5 & 25.6000\\\hline
    \end{tabular}
\end{table}
\begin{table}[ht]\centering
    \caption{Parameter values of the high-traffic ($\rho=0.73$) system with seven
      environment states.
      \label{tab:params_ht}}
    \begin{tabular}{c|r*{7}{c}|}\hline 
      Environment ($z$) & $\la_z$ & $\mu_z$ & $\xi_z$ & $\al_z$
      & $\tht_z$ & $\rho_z$ \\\hline
      1 &   2.0  & 7.0  & 0.50   & 2.0  & 0.05 & 0.4000\\
      2 &   1.2 &  9.9  & 2.01   & 0.2  & 1.50 & 1.2821\\
      3 &   1.7  & 8.5  & 0.05   & 4.1 & 4.20 & 0.2071\\
      4 &   1.2  & 2.7  & 0.30   & 6.8 & 0.40 & 0.5176\\
      5 &   4.6  & 13.1  & 1.90   & 2.1  & 0.10 & 0.7108\\
      6 &   10.2  & 1.1  & 2.70   & 1.5  & 0.90 & 8.2263\\
      7 &   0.3 & 3.9   & 0.10  & 3.2  & 0.50 & 0.1023\\\hline
    \end{tabular}
\end{table}

\pagebreak
Lastly, we define a common random environment with the infinitesimal
generator $Q$ given by
\[\small Q=\begin{bmatrix*}
    -17.5 &   4.5 &   2.6 &   1.1 &   0.0 &   6.1 &   3.2 \\
     5.8  & -32.3 &   3.2 &   7.8 &   4.4 &   8.2 &   2.9 \\
     2.2  &   9.6 & -40.4 &   0.8 &   8.8 &   7.4 &  11.6 \\
     0.1  &   1.7 &   5.1 & -19.8 &   0.0 &  12.9 &   0.0 \\
     6.5  &   0.0 &   8.2 &   8.1 & -27.4 &   3.7 &   0.9 \\
     6.6  &   8.9 &  16.2 &   3.9 &   8.2 & -45.9 &   2.1 \\
     1.8  &   2.8 &   9.5 &   0.8 &   7.9 &   0.0 & -22.8
  \end{bmatrix*}
\]
The Bright and Taylor algorithm is then applied to each of the systems, which are truncated to a maximum orbit size of $R^*=75$.
A representative set of steady-state performance measures for each of the three resulting systems is provided in Table
\ref{tab:performance} for set inventory threshold values of $s=10$ and $S=35$. Within this set of values, it can be verified that the average long-run inventory size $\binv=23$ in Table \ref{tab:performance}, which is computed directly from the
first moment of inventory size Eqn. \eqref{eq:Binv} for each system,
agrees with the value of $\binv$ calculated using formula Eqn. \eqref{eq:newBinv}
of Lemma \ref{lemma:B_inv}.

\addtolength{\tabcolsep}{-4pt}
\begin{table}[ht]
    \centering\small
    \caption{Marginal steady-state probabilities and performance measures, $s=10$, $S=35$.}
    \label{tab:performance}
    \begin{tabular}{lccccccccccc}
  & \multicolumn{3}{c}{\textbf{Probability}} & \multicolumn{6}{c}{\textbf{Performance Measure}}\\ \cline{2-4} \cline{5-11}
\multicolumn{1}{l|}{\textbf{Traffic}}   & \textbf{Idle} & \textbf{Busy}	& \textbf{Failed} &	\multicolumn{1}{|c}{\bm{$L_R$}} & \bm{$L$}	& \bm{$W_R$} & \bm{$W$} & \bm{$\binv$} &  \multicolumn{1}{c|}{\bm{$D_s$}} & \multicolumn{1}{c|}{\bm{$\rho$}} \\ \hline
\multicolumn{1}{|l|}{\textbf{Low}}  & 0.6061 &	0.2052 & 0.1887 & \multicolumn{1}{|c}{1.5971} & 1.8023 & 0.6745 & 0.7612	& 23	&\multicolumn{1}{c|}{19.0299} & \multicolumn{1}{c|}{0.1692}\\ 
\multicolumn{1}{|l|}{\textbf{Medium}} & 0.2577 &	0.5516 & 0.1907 & \multicolumn{1}{|c}{8.8412} &	9.3928 & 3.0707 & 3.2623 &	23 & \multicolumn{1}{c|}{81.5578} & \multicolumn{1}{c|}{0.4071}\\
\multicolumn{1}{|l|}{\textbf{High}} & 0.2891 & 0.4492 & 0.2618 & \multicolumn{1}{|c}{13.3542} &	13.8034 & 4.6129 & 4.7681 &	23 & \multicolumn{1}{c|}{119.2021} & \multicolumn{1}{c|}{0.7305}\\ \hline
    \end{tabular}
\end{table}

\subsection{Results}

The steady-state average cost function will now be defined and then analyzed for the presence of minimal values. The elements of the cost function are defined similarly to those of Ko \cite{Ko2020}, as defined here:
\begin{itemize}
  \item $C_h$: inventory holding cost per item per unit of time,
  \item $C_b$: blocking cost per item sent to the retrial orbit per unit of time,
  \item $C_o$: reordering cost per order from the supplier,
  \item $C_p$: purchase, or procurement, cost per item.
\end{itemize}
Using these elements, together with the average system parameters
defined in Section 5, we define the steady-state mean total cost per unit time $C_T$ as
\begin{align*}
  C_T(s,S) &= C_h\binv + C_bL_R + C_o r_o + C_p r_s \\[1ex]
      &= C_h\apar{\frac{1}{2}(S+s+1)} + C_bL_R + C_o\apar{\frac{1}{(S-s)W}}
         + C_p\apar{\frac{1}{W}} \\[1ex]
      &= \frac{C_h}{2}(S+s+1) + C_b L_R
         + \frac1W \apar{C_p+\frac{C_o}{S-s}},
\end{align*}
and the pairs $\sS$ of inventory thresholds are, for some fixed $a\in\Zp$, taken from the feasible region
\[\mathcal{F}_a=\sett{\sS\in\Zp^2\,|\,s=a,\,S\ge s+1}.\]
The steady-state average cost optimization problem may now be stated  as
\begin{equation}\label{eq:opt_prog}
  \text{Minimize }\; C_T\; \text{ subject to }\; \sS\in\mathcal{F}_a
\end{equation}
where the cost coefficients are assigned the fixed values
\[C_h=5,\quad C_b=24,\quad C_0=11,\quad C_p=2.\]\vspace{-5ex}

By Proposition \ref{prop:ind_of_perf_meas}, it may be inferred
that the performance measures $L_R$ and $W$ are independent
of $s$ and $S$, and are therefore constant in $\mathcal{F}_a$. This permits an analytical solution to optimization problem Eqn. \eqref{eq:opt_prog}, which appears in the following theorem.
\begin{thm}\label{thm:opt_solution}
    The unique optimal solution $\sSstar$ to Eqn. \eqref{eq:opt_prog} is given by
    \begin{equation*}
      \sSstar = \arg \min{ C_T(a, \max \{ a+1, \lfloor \bar{S} \rfloor \} },
        C_T(a, \lceil \bar{S} \rceil) \},
    \end{equation*}
    where
    \[\bar{S} =  a + \sqrt{ \frac{2 C_o}{C_h W}}.\]
\end{thm}
\begin{proof}
     Let $\sS$ be an arbitrary feasible point satisfying $s > a$.  Then for any value of
     $S$, the value of $C_T$ can be reduced by decreasing both $s$ and $S$ by the same amount,
     since all terms in the objective function are held constant except the first term, which is
     reduced.  Therefore, $\sS$ cannot be locally optimal, which means that $s = a$ is optimal.
     This reduces the optimization problem to the single variable minimization of
  \begin{eqnarray*}
     f(S) &=& \frac{C_h}{2} (S + a + 1) + C_b L_R 
           +  \left( C_p + \frac{C_o}{S - a} \right) 
              \left( \frac{1}{W} \right),
  \end{eqnarray*}
  where $S > a$ must hold.
  Then applying first and second-order optimality conditions, we have \vspace{-1ex}
  \begin{eqnarray*}
      f'(S)  &=& \frac{C_h}{2} - \frac{C_0}{W(S-a)^2} = 0, \\
      f''(S) &=& \frac{2 C_0}{W (S-a)^3} > 0,
  \end{eqnarray*}
  for all $S > a$.  Solving for $S$ yields $S = \bar{S}$ and
  $f''(\bar{S}) > 0$.  Since $S$ must be integer and $f''(S) > 0$ for all $S > a$, $f$ is convex
  (which ensures uniqueness in this case),
  and the result is obtained by taking the two integer values that bracket $\bar{S}$ and choosing the one with a smaller function value, but ensuring that the floor function does not drop below $a$ (to enforce $S \geq a + 1$).
\end{proof}

By applying Theorem \ref{thm:opt_solution} to the systems constructed in Section 6.1 for $a=1$, we arrive at the results for an optimal steady-state average cost that appears in Table~\ref{tab:optResults}. The values of the optimal average cost $C_T\sSstar$ demonstrate the expected monotone increasing behavior with traffic intensity, primarily due to the penalty cost $C_b$ for demands held in orbit. Observe that the minimum cost for the low- and medium-traffic systems are not corner points of $\mathcal{F}_a$, in spite of the constant value $s=1$ of the lower inventory threshold. It is anticipated that a similar formulation for an analogous model with replenishment delay and accompanying penalties for such delays will result in interior-point solutions for those models.


\begin{table}
  \centering
  \caption{Optimal $\sSstar$ settings for the low-, medium-, and high-traffic systems, $a=1$}
  \label{tab:optResults}
  \begin{tabular}{l|cc}
     \hline
     Traffic Setting & $\sSstar$ & $C_T\sSstar$ \\ \hline
     Low  &    $(1,5)$      & 42.57      \\
     Medium &    $(1,5)$      & 44.35      \\       
     High   &    $(1,4)$      & 49.49     \\ \hline
  \end{tabular}
\end{table}

\section{Conclusion}

Due to the novel approach enabled by the results contained in Cordeiro et al. \cite{cordeiro2019}, it is now possible to derive a closed-form traffic intensity condition for a complex inventory system with exponential rate parameters modulated by a random environment. Particularly notable is the compact matrix-vector form of the traffic intensity formula, whose complexity of expression is unaffected by the number of defined environments and the magnitude of inventory thresholds, thus enabling the construction and
subsequent numerical investigation of stable systems.

Such follow-on numerical studies must first proceed with the computation of optimal steady-state average costs for systems with replenishment delay. While simplifying the assumption of instantaneous replenishment is sufficient to demonstrate
the efficacy of the method of Cordeiro et al. \cite{cordeiro2019} in deriving a closed-form traffic intensity and to the provision of a basic framework for a cost-optimization study, it prevents the analysis of performance measures that pertain to delays in stock replenishment. It is anticipated that extending the current inventory
model to incorporate such a feature would facilitate a more
comprehensive numerical investigation into its optimal-cost characteristics.

Beyond such incremental directions in the study of inventory systems similar to the one of this paper, the method described herein to derive a closed-form traffic intensity may potentially be used for other queueing inventory models whose underlying Markov chains are row-convergent LDQBDs. Some relevant examples are multi-server queueing models, multiple product inventory models, and perishable systems with Markovian product degradation, among what is anticipated to be many others.


\vspace{3ex}
\noindent\textbf{Competing Interests}:
The authors have no competing interests to report.

\bibliographystyle{abbrv}
\spacingset{1}
\bibliography{retrial_inventory_ldqbd}

\begin{thebibliography}{10}

\bibitem{arrow1951}
K.~J. Arrow, T.~Harris, and J.~Marschak.
\newblock {Optimal Inventory Policy}.
\newblock {\em Econometrica}, 19(3):250--272, Jul 1951.

\bibitem{artalejo_etal06}
J.~Artalejo, A.~Krishnamoorthy, and M.~Lopez-Herrero.
\newblock Numerical analysis of {$(s,S)$} inventory systems with repeated
  attempts.
\newblock {\em Annals of Operations Research}, 141(1):67--83, 2006.

\bibitem{bri_tay95}
L.~Bright and P.~Taylor.
\newblock Calculating the equilibrium distribution in level dependent
  quasi-birth-and-death processes.
\newblock {\em Communications in Statistics: Stochastic Models},
  11(3):497--525, 1995.

\bibitem{cinlar_intro_stoch}
E.~{\c{C}}{\i}nlar.
\newblock {\em Introduction to Stochastic Processes}.
\newblock Prentice-Hall, Englewood Cliffs, NJ, 1975.

\bibitem{cordeiro2019}
J.~Cordeiro, J.~Kharoufeh, and M.~Oxley.
\newblock On the ergodicity of a class of level-dependent quasi-birth-and-death
  processes.
\newblock {\em Advances in Applied Probability}, 51(4):1109--1128, 2019.

\bibitem{Feldman1978}
R.~M. Feldman.
\newblock {Continuous review $(s, S)$ inventory system in a random
  environment}.
\newblock {\em Journal of Applied Probability}, 15(3):654--659, 1978.

\bibitem{fisher_horn_2000}
J.~D.~M. Fisher and A.~Hornstein.
\newblock $({S}, s)$ inventory policies in general equilibrium.
\newblock {\em Review of Economic Studies}, 67(1):117--145, 2000.

\bibitem{IglehartKarlin1963}
D.~Iglehart and S.~Karlin.
\newblock Optimal policy for dynamic inventory process with nonstationary
  stochastic demands.
\newblock In K.~Arrow, S.~Karlin, and H.~Scarf, editors, {\em Studies in
  Applied Probability and Management Science}, pages 127--147. Stanford
  University Press, Redwood City, California, USA, 1962.

\bibitem{Karlin1960}
S.~Karlin.
\newblock Dynamic inventory policy with varying stochastic demands.
\newblock {\em Management Science}, 6(3):231--258, Apr 1960.

\bibitem{Ko2020}
S.-S. Ko.
\newblock {A nonhomogeneous quasi-birth-death process approach for an $(S, s)$
  policy for a perishable inventory system with retrial demands.}
\newblock {\em Journal of Industrial \& Management Optimization},
  16(3):1415--1433, May 2020.

\bibitem{Krishnamoorthy2012}
A.~Krishnamoorthy, S.~Nair, and V.~C. Narayanan.
\newblock {An inventory model with server interruptions and retrials}.
\newblock {\em Operational Research}, 12(2):151--171, Sep 2012.

\bibitem{kulkarni_text}
V.~G. Kulkarni.
\newblock {\em Modeling and Analysis of Stochastic Systems}.
\newblock Chapman \& Hall/CRC Texts in Statistical Science. Taylor \& Francis,
  Boca Raton, FL, 1st edition, 1996.

\bibitem{neuts71}
M.~F. Neuts.
\newblock A queue subject to extraneous phase changes.
\newblock {\em Advances in Applied Probability}, 3:78--119, 1971.

\bibitem{neuts2:78}
M.~F. Neuts.
\newblock Further results on the {$M/M/1$} queue with randomly varying rates.
\newblock {\em OPSEARCH}, 15(4):158--168, 1978.

\bibitem{neuts1:78}
M.~F. Neuts.
\newblock The {$M/M/1$} queue with randomly varying arrival and service rates.
\newblock {\em OPSEARCH}, 15(4):139--157, 1978.

\bibitem{neuts_mgssm_81}
M.~F. Neuts.
\newblock {\em Matrix-Geometric Solutions in Stochastic Models: An Algorithmic
  Approach}.
\newblock Dover Books on Advanced Mathematics. Dover Publications, New York,
  NY, 1981.

\bibitem{Ozekici1999}
S.~{\"{O}}zekici and M.~Parlar.
\newblock Inventory models with unreliable suppliers in a random environment.
\newblock {\em Annals of Operations Research}, 91:123--236, 1999.

\bibitem{Perry2002}
D.~Perry and M.~Posner.
\newblock Production-inventory models with an unreliable facility operating in
  a two-state random environment.
\newblock {\em Probability in the Engineering and Informational Sciences},
  16(3):325--338, 2002.

\bibitem{senn83}
L.~I. Sennott, P.~A. Humblet, and R.~L. Tweedie.
\newblock {Mean Drifts and the Non-Ergodicity of Markov Chains.}
\newblock {\em Operations Research}, 31(4):783--789, 1983.

\bibitem{Song1993}
J.~S. Song and P.~Zipkin.
\newblock {Inventory control in a fluctuating demand environment}.
\newblock {\em Operations Research}, 41(2):351--370, 1993.

\bibitem{ushakumari06}
P.~V. Ushakumari.
\newblock On $(s, {S})$ inventory system with random lead time and repeated
  demands.
\newblock {\em Journal of Applied Mathematics \& Stochastic Analysis}, Volume
  2006. Article ID 81508:1--22, 2006.

\bibitem{yech_naor71}
U.~Yechiali and P.~Naor.
\newblock Queuing problems with heterogeneous arrivals and service.
\newblock {\em Operations Research}, 19(3):722--734, 1971.

\end{thebibliography}
	
\section{Appendix}

{\bfseries\emph{Proof of Lemma \ref{lemma:B_inv}}}
\begin{proof}
  We first define the continuous-time stochastic process
  $\sigma=\sett{I(t)}[t\ge 0]$ on the state-space of inventory states
  $\mcI=\sett{s+1,s+2,\dots,S}$. We will first need to establish that $\sigma$
  is a semi-Markov process (SMP) with transition epochs $S_0, S_1, S_2,\dots$
  taken from the service completion times of $\Phi$, with $S_0=0$.
  To do this, we will define the variables
  \[Y_n=I(S_n+),\quad n\in\Zp.\]
  The process $\sigma$ is considered an SMP if (1) it is a
  piecewise-constant, left-continuous process, and (2) the sequence of bivariate
  random variables $\sett{(Y_n,S_n)}[n\ge 0]$ is a Markov renewal sequence (MRS).
  It can easily be seen that (1) is inherited from $\Phi$. Using
  the well-known fact that end-of-service epochs in a Markovian
  queueing system are stopping times, as defined in \c{C}inlar
  \cite{cinlar_intro_stoch}, (2) may be shown by means of a routine validation of each of the axioms that define an MRS (see Kulkarni \cite{kulkarni_text}).
  We may thus conclude that $\sigma$ is an SMP with kernel $G(x)=[G_{ij}(x)]$,
  where
  \[G_{ij}(x)=\prb{Y_{n+1}=j,\,S_1\le x}[Y_0=i],\quad i,j\in\mcI.\]
  Furthermore, $\sigma$ possesses the embedded DTMC $\eta=\sett{Y_n}[n\ge 0]$ with
  the associated transition probability matrix $P=G(\infty)$.
  
  
  With $\sigma$ established as an SMP, we may now utilize Kulkarni
  \cite[Theorem 9.27]{kulkarni_text} to compute the steady-state distribution
  for $\sigma$. This result requires that $\sigma$ exhibits the properties of irreducibility, aperiodicity, and positive recurrence. Since irreducibility and periodicity are inherited from the parent process $\Phi$, it remains to show that $\sigma$ is positive recurrent.  
  Let $T_j$ be the time of first jump of $\sigma$ to state
  $j\in\mcI$, namely
  \[T_j=\inf\sett{t\ge S_1}[I(t)=j,\,I(t-)\ne j],\]
  which is also the time of the first entry of $\Phi$ into the set of states $C_j\subset S_{\Phi}$ for which $I(t)=j$. Also, define the conditional distributions of time for $\sigma$ to reach state $j$ from $i$ and the expectations associated with these distributions as
  \begin{align*}
    F_{ij}(t) &= \prb{T_j\le t}[I(0)=i], \\
    \mu_{ij}  &= \E{T_j}[I(0)=i].
  \end{align*}
  These may likewise be interpreted as the conditional probabilities of the time of the first jump of $\Phi$ into $C_j$ beginning in $C_i$, and their expected values. In order to conclude that $\sigma$ is a positive recurrent SMP, it must be shown that $F_{ii}(\infty)=1$ and $\mu_{ii}<\infty$ for every $i\in\mcI$.
  
  Consider any state $i\in\mcI$ and suppose that $\sigma$ is
  in state $i$ at time $t=S_0=0$. In this case, $\Phi$ is presumed to be in some state
  \[y=Y(0)=(R(0),i,X(0),Z(0))\in C_i\subset\Sphi\]
  where $Y(t)=(R(t),I(t),X(t),Z(t))$.
  After the return time $T_i$ has elapsed, $I(T_i)=i$ for the SMP $\sigma$.
  However, it is possible that $\Phi$ is in a different state
  \[y'=Y(T_i)=(R(T_i),i,X(T_i),Z(T_i))\in C_i\subset\Sphi,\quad y\ne y',\]
  and thus $T_i\le T^{\Phi}_{y}$, which is the time of first return of $\Phi$ to $y\in C_i$. However, if we consider the quantities
  \begin{align*}
    F^{\Phi}_{yy'}(t) &= \prb{T^{\Phi}_{y'}\le t}[Y(0)=y], \\
    \mu^{\Phi}_{yy'}  &= \E{T^{\Phi}_{y'}}[Y(0)=y],
  \end{align*}
  then, since $\Phi$ has been assumed to be positive recurrent, it
  must then be true that $F^{\Phi}_{yy}(\infty)=1$ and $\mu^{\Phi}_{yy}<
  \infty$. In other words, $\Phi$ returns to state $y\in C_i$ with probability
  1, which simultaneously implies that $\sigma$ must likewise return to
  $i$ with probability 1. Hence $F_{ii}(\infty)=1$. Moreover, since
  $T_i\le T^{\Phi}_y$ for any given initial state
  $y\in C_i$ for $\Phi$, the properties of expected values yield the inequality
  \[\mu_{ii}\le\mu^{\Phi}_{yy}<\infty.\]
  Therefore, since $i\in\mcI$ was arbitrarily chosen, $\sigma$ must be positive
  recurrent.
  
  It now remains to compute the steady-state probability
  distribution of the SMP $\sigma$, from which we may obtain the quantity
  $\binv$. For each $j\in\mcI$, let
  \[p_j=\lim_{t\to\infty}\sett{I(t)=j}[I(0)=i]\]
  be the steady-state probability of being in state $j$ for the irreducible,
  aperiodic, and positive recurrent SMP $\sigma$, which is then given by the
  expression
  \[p_j=\frac{\pi_j\mu_j}{\sum_{k=s+1}^S\pi_k\mu_k},\]
  where $\bpi=[\pi_j]$ is a positive row vector solution to the system
  $\bpi P=\bpi$, if it exists, and $\mu_j$ is the expected sojourn time of $\sigma$
  in state $j\in\mcI$.
  To compute $\bpi$, we first construct the matrix $P$ of the embedded DTMC $\eta$, which becomes
  \[P=G(\infty)=
    \begin{bmatrix*}
      0 & 0 & \dots & 0 & 0 & 1 \\
      1 & 0 & \dots & 0 & 0 & 0 \\
      0 & 1 & \ddots & 0 & 0 & 0 \\
      \vdots & \vdots & \ddots & \ddots & \vdots \\
      0 & 0 & \dots & 1 & 0 & 0 \\
      0 & 0 & \dots & 0 & 1 & 0 \\
    \end{bmatrix*}.
  \]
  It is possible to visually determine that $\bpi=\eins$; 
  that is, $\pi_j=1$ for each $j\in\mcI$.
  
  In order to determine the quantities $\mu_j$, we observe that
  the i.i.d. successive durations of time between service completions
  $S_{n+1}-S_n$ for each $n\ge 0$ coincide with the sojourn times of $\sigma$ in each
  of its states $j\in\mcI$. Moreover, since Proposition \ref{prop:ind_of_perf_meas}
  informs us that the length of these sojourn times is independent of
  any given inventory size $j$, we may conclude that
  \[\mu_j=\tau_s=\E{S_{n+1}-S_n}[n\ge 0]=\E{S_1-S_0}\text{ for each } j\in\mcI.\]
  Also, because of the positive recurrence of $\sigma$, we have
  \[\tau_s\le\mu_{jj'}\le\mu_{jj}<\infty\text{ for each } j,j'\in\mcI.\]
  Therefore, the steady-state probability of inventory size
  $j\in\mcI$ may be calculated as
  \[p_j=\frac{\tau_s}{\sum_{k=s+1}^S\tau_s}=\frac{1}{S-s}.\]
  Substituting each of these terms into Eqn. \eqref{eq:Binv} gives
  \[\binv=\frac{1}{S-s}\sum_{I=s+1}^SI,\]
  whereupon application of the identity
  \[1+2+\dots+N = \frac{1}{2}N(N+1),\quad N=1,2,3,\dots,\]
  yields Eqn. \eqref{eq:newBinv}.
\end{proof}

\end{document}